\documentclass[12pt]{article}
\usepackage{latexsym,amsmath,amssymb}
\usepackage{epsfig,graphics}
\advance\topmargin -1cm
\advance\oddsidemargin -1.5cm
\newcommand{\cqd}{\hfill\rule{2mm}{2mm}}
\def\car{\mbox{\rm{1\hspace{-0.10 cm }I}}}

\newcommand{\st}{^{(t)}}

\newcommand{\stM}{^{(t+1)}}
\newtheorem{prop}{Proposition}
\newtheorem{teor}{Theorem}

\newtheorem{lema}{Lemma}
\makeatletter

\@addtoreset{equation}{section}
\makeatother
\author{Ricardo R\'{\i}os\footnote{Universidad Central de Venezuela, Facultad de Ciencias,
Escuela de Matem\'aticas, Caracas 1040, Venezuela. Email: rrios@euler.ciens.ucv.ve.}
and
Luis Rodr\'\i guez\footnote{Universidad de Carabobo, Facultad de Ciencias y Tecnolog\'{\i}a,
Departamento de Matem\'aticas, Valencia, Venezuela. Email: larodri@uc.edu.ve.}
}
\title{Estimation in autoregressive models
with Markov regime}
\begin{document}
\maketitle

\begin{abstract}
In this paper we derive the consistency of the penalized
likelihood method for the number state of the hidden Markov chain
in  autoregressive models with Markov regimen. Using a SAEM type
algorithm to estimate the models parameters. We test the null
hypothesis of hidden Markov Model against an autoregressive
process with Markov regime.

\noindent\textbf{Keywords:} Autoregressive process, hidden Markov,
switching, SAEM algorithm, penalized likelihood.
\end{abstract}

\section{Introduction}
This paper is devoted to estimate of autoregressive models with
Markov regime. Our goals in this paper are:
\begin{itemize}
    \item   Estimate, using maximum likelihood estimation (MLE)
        methods, the parameters that define the functions, the
        transition probabilities of the
        hidden Markov chain and the noise variance, computed via SAEM,
        a stochastic version of EM algorithm \cite{Lavielle}, for
        a pre-fixed  number states of the hidden Markov chain.
    \item Test the null hypothesis of HMMs against AR-RM.
    \item Derive the consistency of the penalized likelihood method
    for the number of state.
\end{itemize}

An autoregressive model with Markov regime
(AR-MR) is a discrete-time process defined by:
\begin{equation}
    \label{model}
    Y_n=f_{X_n}(Y_{n-1},\theta_{X_n})+\sigma\varepsilon_n
\end{equation}
where $\{X_n\}$ is a Markov chain with finite state space  $ \{1,
\ldots, m\}$. The transition probabilities denoted as $a_{ij}
=\mathbb{P}(X_n=j|X_{n-1}=i)$. The $a_{ij}$  form  an $m\times m$
transition matrix A. The functions $\{f_1,\ldots,f_m\}$ belong to
a parameterized family
\begin{equation}
\{\theta_1(\cdot)+\theta_0: \theta=(\theta_1,\theta_0)\in\Theta\}
\end{equation}
where  $\Theta$ a compact subset of $\mathbb{R}^2$, and
$\{\varepsilon_n\}$ is a sequence of independent identically distributed standard
normal random variables, $\mathcal{N}(0,1)$. As process
$\{X_n\}$ is not observable then we are forced to work
with simulations of the law of the hidden chain  and to rely on
observed data $\{Y_n\}$ for any inference task.

The usage of Markov regime offers possibilities for modelling time
series ``{\it subject to discrete shifts in regime-episodes across
which the dynamic behavior of the series is markedly different}",
as noted by Hamilton \cite{Hamilton} who used  a model AR-RM in
the context of econometrics,
 for the analysis of
 the U.S. annual GNP (gross national product) series, with two regimes:
contraction and expansion. Linear autoregressive process with Markov regime
are also used in several electrical engineering areas including tracking of
manoeuvring targets, failure detection and stochastic adaptative control
(Douc {\it et alii} \cite{douc}).

An important class of AR-MR is the hidden Markov models (HMMs) for
which the
 functions $\{f_1,\ldots,f_m\}$
are constants ($\theta_{1,i}=0$, for all $i\in\{1,\ldots,m\}$). The HMMs are
used in many different areas: basic and applied sciences, industry,
economics, finance, images reconstruction, speech recognition,
tomography, inverse problem, etc. \cite{Cappe}, \cite{Zucchini}.

The advantage of using the SAEM algorithm  is easiness of movement
in different modal areas, that reduces the chance of the estimate
to avoid a local maximum. The particularities of our problem
allows us to do an exact simulation  of the distribution of the
hidden chain conditional to the observations, using  Carter-Kohn
algorithm \cite{Carterkohn}.

For the hypothesis test of HMM against Linear AR-RM we follow the
ideas of Giudici {\it{et al}} \cite{giudici} then we obtain the
usual asymptotical theory. They used likelihood-ratio test for
HMMs,  to establish that the standard asymptotic theory rests
valid. They work with hidden graphical Gaussian models.

When the number $m$ is unknown, the hypothesis test with
likelihood ratio  techniques fails to estimate $m$ because
regularity hypothesis are not satisfied. Particulary, the model is
not identifiable, in the sense of Dachuna-Duflo \cite{dacuna}
(227), so standard  $\chi^2$  can not be applied.

In the HMM framework, we distinguish two cases according if the
number state of the observed variables is finite or not. In the
finite case, Finesso \cite{finesso} gives $\hat{m}$ a strong
consistence penalized estimator of $m$, assuming that $m$ belongs
to a bounded subset of the integers numbers. Liu and Narayan
\cite{Liu}, also assume this bounded condition and postules a
strongly consistent and efficient  $\hat{m}$ with the probability
of underestimation decaying exponentially fast w.r.t. $N$, while
the probability of overestimation does not exceed $cN^3$. Gassiat
and Boucheron \cite{gassiat3} prove the strong consistence of a
penalized $\hat{m}$ without assumptions about upper bounds for
$m$, with the probability of underestimation and overestimation
decaying exponentially fast.
 In the non-finite case, the likelihood
ratio is not bounded, Gassiat and  Keribin studies in
\cite{gassiat1} show divergence to infinity. As far as we know,
the divergence rate rests unknown. In Gassiat \cite{gassiat2}
results over penalized likelihood are given in order to obtain
weak consistence for the estimator of the number state. We obtain
strong consistence for a penalized $\hat{m}$ in a linear AR-MR,
and $m$ in a bounded set.

The paper is organized as follows. Main assumptions are given in Section 2.
In Section 3 for a fixed number state of the hidden Markov chain,
an SAEM type algorithm is used to estimate
the parameters and is  present the method of simulation of
the hidden Markov chain
and their convergence properties.
In the Section 4 we presents our results on the analysis of LR test.
In
Section 5 we derive the consistency of the penalized likelihood
method for a number state problem. For sake the clarity the
proof of the Lemma 1 is relegated to Appendix A.
Appendix B is devoted to simulations.

\section{Notation and assumptions}
Let $y_{1:N}=(y_1,\ldots,y_N)$ denote the observations and
$X_{1:N}=(X_1,\ldots,X_N)$ the associated vector of the hidden
variables. Using $p$ as generic symbol for densities and distributions,
the likelihood function is given  by
\begin{equation}
p(y_{1:N}|y_0,\psi)=\sum_{x\in\{1,\ldots,m\}^N}p(y_{1:N},X_{1:N}=x|y_0,\psi),
\label{vero1}
\end{equation}
where $x=(x_1,\ldots,x_N)$ and  $\psi=$
$(A,\theta,\sigma)\in\Psi$ and $\Psi =[0,1]^{m^2}\times\Theta^m\times(\mathbb{R}^+)$.
The maximum likelihood estimate (MLE) $\hat{\psi}$ is defined
as,
$$
 \hat{\psi}=\underset{\psi\in\Psi}{\mbox{argmax}}\ \!p(y_{1:N}|y_0,\psi).
$$

Suppose that $Y_0$, $\{X_n\}$ and $\varepsilon_1$ are mutually independent then
\begin{equation}
  p(y_n | x_n, \ldots, x_0, y_{n-1}, \ldots, ,
  y_{0}) = p(y_n | x_n, y_{n-1}).
  \label{indcondicional}
\end{equation}
Using (\ref{indcondicional}) and from Markov property of $\{X_{n}\}$ we have
\begin{align}
     \nonumber
    {l_N(\psi)} &= \log p(y_{1:N}|y_0,\psi)\\
      \nonumber
      &= \log\left( \sum_{x\in \{1,\ldots,m\}^N}{p}(y_{1:N},X_{1:N}=x_{1:N}|\psi) \right)\\
      \nonumber
      &= \log\left( \sum_{x_{1:N} \in \{1,\ldots,m\}^N}
      p(y_{1:N}| X_{1:N}=x_{1:N}, y_0, \psi) p(X_{1:N}=x_{1:N}|\psi)\right)\\
      &= \log\left( \sum_{x_1=1}^m\ldots\sum_{x_m=1}^m
      \prod_{n=1}^Np(y_n|X_{n}=x_n,y_{n-1})\prod_{n=1}^{N-1}a_{x_{n},x_{n+1}}p(X_1=x_1)\right)
      \label{vero}
\end{align}
with
    \begin{eqnarray*}
    p(y_{n}|y_{n-1},X_{n}=i)
     =\frac{1}{\sqrt{2\pi\sigma^2}}
    \exp\left(-\frac{(y_{n}-f_{i}(y_{n-1},\theta_{i}))^2}{2\sigma^2}\right).
    \end{eqnarray*}

For the consistence the MLE we will assume the followings conditions,
\begin{enumerate}
\item[(C1)] The transition probability $A$ is positive, this is,
$a_{ij}\geq \delta$, for all $i,j\in\{1,\ldots,m\}$ for some $\delta>0$.
\end{enumerate}
This condition implies that there is an unique invariant distribution
  $\mu=(\mu_1,\ldots,\mu_m)$.
\begin{enumerate}
\item[(C2)] Let  $
\sum_{i=1}^m\log|\theta_{1,i}|\mu_i<0.$
\end{enumerate}
This condition, and the existence of the moments the
$\varepsilon_1$, implies that the chain extended $\{(Y_n,X_n)\}$
is a geometrically ergodic Markov chain on the  state space
$\mathbb{R}\times\{1,\ldots,m\}$ under $\psi_0$ (see Yao and
Attali \cite{Yao}).
\begin{enumerate}
\item[(C3)]Let $b_{+}:=\sup_{\psi}\sup_{y_0,y_1,i}p(y_1|y_0, i)<\infty$ and
$\mathbb{E}(|\log b_{-}(y_1,y_0)|)<\infty$, where $b_{-}(y_1,y_0):=
\inf_{\psi}\sum_{i=1}^m p(y_1|y_0,i)$.
\item[(C4)]For all $i,j\in\{1,\ldots,m\}$ and all $y,y'\in\mathbb{R}$,
the functions $\psi\to a_{ij}$ and $\psi\to p(y'|y,i)$ are continuous.
\item[(C5)]The model is identifiable in the sense that
$p_\psi=p_{\psi*}$  implies that $\psi=\psi*$.
For this is sufficient that $\theta_i\not=\theta_j$ if $i\not=j$,  up to
an index  permutation (Krisnamurthy and Yin \cite{Kris}).
\item[(C6)] For all $i,j\in\{1,\ldots,m\}$ and $y,y'\in\mathbb{R}$,
the functions $\psi\to a_{ij}$
and $\psi \to p(y'|y,i)$ are twice continuously differentiable over
$O=\{\psi\in\Psi:\ |\psi-\psi_0|<\delta\}$.
\item[(C7)] Let us denote $\nabla$ for gradient operator and
$\nabla^2$ for Hessian matrix,
    \begin{enumerate}
    \item $\sup_{\psi\in O}\sup_{i,j}\|\nabla\log a_{ij} \|<\infty$ and
    $\sup_{\psi}\sup_{i,j}\|\nabla^2\log a_{ij}\|<\infty$.
    \item $\mathbb{E}_{\psi_0}\left(\sup_{\psi\in O}\sup_{i,j}\|
    \nabla\log p(Y_1|Y_0,i)\|^2\right)<\infty$ and
    $\mathbb{E}_{\psi_0}
    \left(\sup_{\psi\in O}\sup_{i,j}\| \nabla^2\log p(Y_1|Y_0,i)\|\right)$
    $<\infty$.
    \end{enumerate}
\item[(C8)]
    \begin{enumerate}
    \item For all $y,y'\in\mathbb{R}$ there exist an integrable  function
    $h_{y,y'}:\{1,\ldots,m\}\to\mathbb{R}^+$  such that
    $\sup_{\psi\in O}p(y_1|y_0,i)\leq h_{y,y'}(i)$.
    \item For all $y,y'\in\mathbb{R}$ there exists  integrable functions
    $h_{i,y}^1:\mathbb{R}\to\mathbb{R}^+$ and $h_{i,y}^2:\mathbb{R}\to\mathbb{R}^+$
    such that $\|\nabla\log p(y'|y,i)\|\leq h_{i,y}(y')$ and
    $\|\nabla^2\log p(y'|y,i)\|\leq h_{x,y}(y')$ for all $\psi\in O$.
    \end{enumerate}
\end{enumerate}
In the next proposition we collect some the results of Douc {\it
et alii} \cite{douc} that attains our work.

\begin{prop}
\label{douc}
\mbox{}
\begin{enumerate}
\item[i)] Assuming (C1)-(C4). Then
$$\lim_{N\to\infty}\sup_{\psi\in\Psi}\left|N^{-1}
l(\psi) -H(\psi)\right|,\ \mathbb{P}_{\psi_0}-a.s$$
where $H(\psi)=\mathbb{E}_{\psi_0}(\log p(Y_0|Y_{-\infty:-1},\psi_0))$.
\item[ii)] Assuming (C1)-C5). Then
$$
\lim_{N\to\infty}\hat{\psi}_{N}=\psi_0\ \mathbb{P}_{\psi_0}-a.s,
$$
\item[iii)] Assuming (C1)-(C3) and (C6)-(C8) then,
$$
N^{1/2}\nabla^2_\psi l(\psi)\to I(\psi_0)\ \mathbb{P}_{\psi_0}-a.s.
$$
\item[iv)] Assuming (C1)-(C8) and that the Fisher information matrix for $\{Y_n\}$,  $I(\psi_0)$
 is positive definite. Then
$$
N^{1/2}(\hat{\psi}_{N}-\psi_0)\to\mathcal{N}(0,I(\psi_0)^{-1})\
\mathbb{P}_{\psi_0}-weakly.
$$
\end{enumerate}
\end{prop}

\section{The estimation algorithm for fixed $m$}
Since the likelihood estimator $\hat{\psi}$ is a solution the
equation $\nabla_\psi l(\psi)=0$,
and this equation do not has an  analytic solution, then
 the maximization has to be performed
numerically by considering  $m^N$ terms in the equation
(\ref{vero}). This restricts the model to observations with
limited size and few states. For HMMs models in a finite space
state  Baum {\it et alii} \cite{Baum70} introduced a
forward-backward  algorithm as an early version of the EM
algorithm.
 The EM algorithm  was proposed by Dempster {\it et alii}
\cite{dempster} to maximize  $\log$-likelihood with
missing data. It enables, with a recursive method,
to change the problem of maximizing the log-likelihood into the problem of maximizing
 some  functional of the
completed  the likelihood $p(y_{1:N},x_{1:N}|\psi)$
of the model:
    $$
    \prod_{n=1}^N\left[\prod_{i,j=1}^m
    a_{ij}^{\car_{i,j}(X_n,X_{n-1})} \prod_{i=1}^m
    p(y_n|y_{n+1},i)^{\car_{j}(X_n)}
    \right],
    $$
where $\car_A(\cdot)$ denotes the function indicator over the set
$A$ and $\car_{A\times
B}(\cdot,\cdot)=\car_A(\cdot)\car_B(\cdot)$. Let us describe  the
$t+1$ step  of the algorithm. Set
$$
Q(\psi,\psi^{(t)})=
\mathbb{E}(\log p(Y_{1:N},X_{1:N}|\psi^{(t)})|Y_{1:N}=y_{1:N},\psi).
$$

Then $Q$ is the expectation of the log-likelihood of the complete data conditioned
to the observed data and the value of the parameter computed at the step $t$,
 $\psi\st$. Then we have that  $Q(\psi,\psi^{(t)})$ equals to
\begin{eqnarray}
\label{Q}
\nonumber
    &&\sum_{n=1}^{N-1}
    \sum_{i,j=1}^m\mathbb{E}(\car_{i,j}(X_n,X_{n+1})|Y_{1:N}=y_{1:N},\psi)\log(a_{ij})\\
    &&\hspace{4cm}+\sum_{n=1}^{N-1}\sum_{i=1}^m\!\mathbb{E}(\car_{i}(X_n)|Y_{1:N}=y_{1:N},\psi) \log
    p(y_{n+1}|y_{n},i).
\end{eqnarray}

The EM is a two steps algorithm: the E step and the M step.
In the E stage compute $Q(\psi,\psi^{(t)})$ the
expectation conditioned to the observed data and the current value of the parameter.

In the M step choose,
$$
\psi^{(t+1)}=\underset{\psi\in\Psi}{\mbox{argmax}}\ \!Q(\psi,\psi^{(t)}).
$$
The EM algorithm converges to a maximum-likelihood estimate for
any initial value, when the complete data likelihood function is
in the exponencial family  and a differentiability condition is
satisfied.

In order to avoid local minima, we have used an stochastic approximation of the EM algorithm,
the SAEM algorithm. Such algorithm has been developed
by Celeux {\it et alii}. in \cite{Celeux83}, \cite{Celeux92} and
\cite{Celeux95}, and its convergence has been proved by
Del\-yon {\it et alii} \cite{Lavielle}.
The EM algorithm is modified in the following way:
the (E) step is split into a simulation step (ES)
and stochastic approximation step (EA):
\begin{enumerate}
\item[ES] Sample one realization $x_{1:N}^{(t)}$ of the missing
data vector under $p(x_{1:N}|y_{1:N},\psi^{t-1})$. \item[EA]
Update the current approximation of the EM intermediate quantity
according to:
    $$Q_{t+1}=Q_t+\gamma_t\left(\log
    p(y_{1:N},x_{1:N}^{(t)}|\psi')-Q_t\right)
    $$
\end{enumerate}
where $(\gamma_t)$ satisfies the condition:
\begin{enumerate}
\item[(RM)]
for all $t\in\mathbb{N},$ $\gamma_t\in[0,1]$, $\sum_{t=1}^\infty\gamma_t=\infty$ and
$\sum_{t=1}^\infty\gamma_t^2<\infty$.
\end{enumerate}

\subsection{ES step}
\label{pasoES}
In this section we describe the simulating method used in
the SAEM algorithm. For sampling under
the  conditional distribution,
    $$p(x_{1:N}|y_{1:N},\psi)=\mu_{x_1}
    p(y_1|y_0,x_1)
    \ldots
    a_{x_{N-1}x_{N}}p(y_N|y_{N-1},x_N)/ p(y_{1:N}|\psi),
    $$
for any  $x_{1:N}=(x_1,\ldots,x_N)\in\{1,\ldots,m\}^N$,  Carter and Kohn
in \cite{Carterkohn} give a method that constitutes a stochastic version
of the forward-backward algorithm
proposed by  Baum {\it et alii} \cite{Baum70}.
This follows by observing that
$p(x_{1:N}|y_{1:N},\psi)$
can be decomposed as,
    $$p(x_{1:N}|y_{1:N},\psi)=
    p(x_N|y_{1:N},\psi)
    \prod_{n=1}^{N-1}p(x_n|x_{n+1},y_{1:N},\psi).
    $$
Provided that $X_{n+1}$ is known, $p(X_n|X_{n+1},y_{1:N},\psi)$
is a discrete distribution, which suggests the following sampling strategy. For $n = 2,
\ldots, N$, $i\in\{1,\ldots,m\}$, compute and store recursively the optimal filter
$p(X_n|y_{1:n},\psi)$ as
    $$
    p(X_n=i|{y}_{1:n},\psi)
    \propto p(y_n|y_{n-1},X_n=i,\psi)
    \sum_{i=1}^{m}a_{ij}p(X_{n-1}=j|{y}_{1:{n-1}}).
    $$

Then, sample $X_N$ from $p(X_N|y_{1:N},\psi)$ and for
$n=N-1,\ldots,n$, $X_n$ is sample from
    $$
    p(X_n=i|X_{n+1}=x_{n+1},{y}_{1:n},\psi)
    =\displaystyle\frac
    {a_{ij_{n+1}}p(X_n=i|{y}_{1:n},\psi)}
    {\sum_{l=1}^{m}a_{il}p(X_n=l|{y}_{1:n},\psi)}.
    $$

As a consequence, the estimation procedure generate an ergodic Markov chain
$\{x_{1:N}^{(t)}\}$ on the finite state space $\{1,\ldots,m\}^N$, so that
$p(x_{1:N}|y_{1:N},\psi)$ is its stationary distribution. Ergodicity
follow from irreducibility and aperiodicity, by observing the positivity of the
kernel, this is,
    \begin{eqnarray*}
    K(x_{1:N}^{(t)}|x_{1:N}^{(t-1)},\psi)\propto&p(x_N^{(t)}|\psi,y_{1:N})
    \prod_{n=1}^{N-1}
    p(x_n^{(t)}|x_{n+1}^{(t)}\psi,y_{1:N})>0.\
    \end{eqnarray*}

In this case the standard ergodic result for finite Markov chains applies
(for instance, Kemeny and
Snell \cite{Kemeny}),
$$
\|K(x_{1:N}^{(t+1)},x_{1:N}^{(t)},\psi)-p(X_{1:N}|y_{1:N},\psi)\|\leq C\rho^{t-1},
$$
with $C=card(\{1,\ldots,m\}^N)$, $\rho=(1-2K_x^*)$ y $K^*=\inf
K(x'|x,\psi)$, for $x,x'\in\{1,\ldots,m\}^N$.

\subsection{EA step}
The (\ref{Q}) equation suggests us to substitute the step EA for approximations of
Robins Monro (ver Duflo \cite{duflo}),
$s=(s_{1}\stM,s_{2}\stM,s_3\stM)$,
defined by:
    \begin{eqnarray}
    \label{estadistico1}
    s_{1}\stM(i,n)&=&s_{1}\st(i,n)+\gamma_t\left(\car_{i}(x_n)-s_{1}\st(i,n)\right)\\
    \label{estadistico2}
    s_{2}\stM(i)&=&s_{2}\st(i)+\gamma_t\left(N_i
    (x_{1:N})-s_{2}\st(i)\right)\\
    \label{estadistico3}
    s_{3}\stM(i,j)&=&s_{3}\st(i,j)+\gamma_t\left(N_{ij}
    (x_{1:N})-s_{3}\st(i,j)\right).\
    \end{eqnarray}
where
$N_i(x_{1:N})=\sum_{n=1}^{N-1}\car_i(x_n)$ and $N_{ij}(x_{1:N})
=\sum_{n=1}^{N-1}\car_{i,j}(x_n,x_{n+1})$,
are sufficient statistics  for the chain of hidden Markov.

When $f_{j}(y,\theta_{j})=\theta_{j}$,
the maximization step is given by,
\begin{eqnarray*}
    \displaystyle\hat{a}_{ij}\stM&=&\frac{s_3\stM(i,j)}{s_2\stM(i)}\\
      \hat{\theta}_i^{(t+1)}
    &=&\frac{\sum_{n=1}^Ns_1\stM(i,n)y_n}{s_2\stM(i)}\\
        \hat{\sigma^2}^{(t+1)}
    &=&\frac{1}{N-1}\sum_{n=1}^{N-1}s_1\stM(i,n)(y_n-\hat{\theta}_i\stM)^2,
\end{eqnarray*}
and for $f_{j}(y,\theta_{j})=\theta_{1,j}y+\theta_{0,j}$ by,
\begin{eqnarray*}
    \displaystyle\hat{a}_{ij}\stM&=&\frac{s_3\stM(i,j)}{s_2\stM(i)}\\
    \hat{\theta}_{1,i}\stM&=&
    \frac{\sum_{n=1}^{N-1}s_1\stM(i,n)y_ny_{n-1}-\sum_{n=1}^{N-1}s_1\stM(i,n)y_n
    \sum_{n=1}^{N}s_1\stM(i,n)y_{n-1}}
    {\sum_{n=1}^{N-1}s_1\stM(i,n)y_{n-1}^2-\left(\sum_{n=1}^{N-1}s_1\stM(i,n)y_n\right)^2}\\
    \hat{\theta}_{0,i}\stM&=&
    \sum_{n=1}^{N-1}s_1\stM(i,n)y_n-\hat{\theta}_{1,i}\sum_{n=1}^{N}s_1\stM(i,n)y_{n-1}\\
    \hat{\sigma^2}^{(t+1)}
    &=&\frac{1}{N-1}\sum_{n=1}^{N-1}\sum_{i=1}^ms_1\stM(i,n)(y_n-f_i(y_{n-1}.
    \hat{\theta}_i))^2
\end{eqnarray*}

 We consider the observations $y_{1:N}$ fixed,
 the previous expressions define, in an explicit way,
 in each one of the two cases of study, the application
 $\hat{\psi}=\psi(s)$ between the sufficient statistics and the parameters space
 necessary to SAEM.

\subsection{Convergence}
The simulation procedure generates
$\{x_{1:N}\st\}$, a finite Markov chain. The
hypotheses of Delyon {\it et alii}
\cite{Lavielle} that ensures the convergence of the SAEM
algorithm
are no more satisfied but
in this case, we can be use the Theorem 1 of
Kuhn and Lavielle in \cite{KuhnLavielle}:
\begin{teor}
 If we suppose the conditions that guarantee the convergence of the
 EM algorithm, the condition (RM) and the following hypothesis,
\begin{enumerate}
\item[SAEM1]The function $p(y_{1:N}|\psi)$ and the function  $\hat{\psi}=\psi(s)$
 are $l$ time differentiable.
\item[SAEM2] The function $\psi\to K_\psi=K(\cdot|\cdot,\psi)$
is continuously differentiable on $\Psi$.
The transition probability  $K_\psi$ generates
a geometrically ergodic chain with invariant probability
$p(x_{1:N}|y_{1:N},\psi)$. The chain $\{x_{1:N}\st\}$
takes  values a compact subset.
The function $s$ is bounded.
\end{enumerate}
Then, w.p 1,
 $lim_{t\to\infty}d(\psi\st,\mathcal{L})=0$ where
$\mathcal{L}=\{\psi\in\Psi:\ \nabla_\psi l(\psi)=0\}$ is the set
of stationary points.
\end{teor}

In our case the hypotheses of the theorem are verified, in fact,
the hypothesis RM is satisfied choosing the sequence $\gamma_t=1/t$,
SAEM1 is obtained because $\varepsilon_1$ is distributed normal and
SAEM2 is consequence of the discussion
made in \S\ref{pasoES}. This guarantees the previous theorem and this give us the
convergency.

\section{Hypothesis test}
In this section we study the likelihood ratio test (LRT) for testing a model HMMs
 against a process AR-RM. We prove that
the standard theory for LRT of a point null hypothesis is valid.
Let $\psi=(A,\theta_1,\theta_0,\sigma^2)$ and
$\psi_0=(A,0,\theta_0,\sigma^2)$, then the test we consider is
that
$$
H_0:\ \theta_1=0
$$
against
$$
H_1:\ \theta_1\not=0.
$$

\begin{teor}
Assume that (C1)-(C8) hold. Then,
$$
2(l(\hat{\psi})-l(\psi_0))\to\chi^2_{m},
$$
under $\mathbb{P}_{\psi_{0}}$.
\end{teor}

\noindent{\bf Proof:} Using the Taylor
expansion of  $l(\psi)$ around $\hat{\psi}$,
$$
l(\psi_0) - l(\hat{\psi}) = (\psi_0-\hat{\psi})\nabla_\psi l(\hat{\psi})
+\frac{1}{2}(\psi_0-\hat{\psi})^t\nabla^2_\psi l(\tilde{\psi})(\psi_0-\hat{\psi})
$$
where $\tilde{\psi}=\lambda\psi_0+(1-\lambda)\hat{\psi},$ $\lambda\in (0,1)$. Also
$\nabla_\psi l(\hat{\psi})=0$. So
$$
-2(l(\psi_0) - l(\hat{\psi}))=
-[N^{1/2}(\psi_0-\hat{\psi})^t][N^{-1}\nabla^2_\psi l(\tilde{\psi})]
[N^{1/2}(\psi_0-\hat{\psi})]
$$
Now, since $\hat{\psi}_N\to\psi_0$ $\mathbb{P}_0-a.s.$ does $\tilde{\psi}_N$,
and using Proposition \ref{douc}-(iii-iv),
$$
N^{1/2}(\hat{\psi}_N-\psi_0)\to\mathcal{N}(0,I(\psi_0)^{-1})\
\mathbb{P}_{\psi_0}-weakly
$$
and
$$
-N^{-1/2}\nabla^2_{\psi}l(\tilde{\psi})\to I(\psi_0)\
\mathbb{P}_{\psi_0}-a.s.
$$

So the proof is complete.\cqd

The theorem says that we can employ the LRT test rejects $H_0$ if:
$$
-2(l(\psi_0) - l(\hat{\psi}))\geq\chi^2_{m,\alpha}
$$
where $\chi^2_{m,\alpha}$ is the $\alpha$-quantile of the $\chi^2_m$ distribution.
\section{Penalized estimation of the number state}
In this section we presents a penalized likelihood method for
selecting  the number state $m$ of the hidden Markov chain
$\{X_n\}$. For each value of $m\geq 1$, we consider the sets
$\Psi_m$ and  $\mathcal{M}=\bigcup_{m\geq 1}\Psi_m$, the
collection of all the different models. For a fixed $m$, we have
seen in Section 3 that it is possible to estimate the unknown
parameters for the model. Hence, it is now possible evaluate the
log-likelihood chosen model $l(\hat{\psi}_m)$.

As we assumed identifiability  $(C5)$, we have that true number state,
 $m_0$ is minimal, that is,
there does not exist a parameter $\psi_{m}\in\Psi_{m}$ with
$m<m_0$ such that $\psi_m$ and $\psi_{m_0}$ induce an identical
law for $\{Y_n\}_{n\geq0}$. We said that $\hat{m}_n$
over-estimate the number state $m_0$ if $\hat{m}_N>m_0$ and
under-estimate the number state if $\hat{m}_N< m_0$.

The penalized maximum likelihood (PML) is defined as:
$$
   C(N,m)=-\log p(y_{1:N}|y_0,x_1\hat{\psi}(N))+pen(N,m),
$$
where $\hat{\psi}(N)$ is the maximum likelihood of $\psi\in\Psi_m$
based on $N$ observations and
$pen(N,m)$ is a positive and increasing function of $m$. A number state estimation
procedure is defined as follows:
$$
\hat{m}(N)=\min\{\underset{m\geq 1}{\mbox{argmin}}\ \!C(N,m)\}.
$$

In the following theorem we prove that the
estimator PML over-estimate the number state $m_0$.

\begin{teor}
Assume (C1)-(C5) and that $\lim_{N\to\infty}pen(N,m)=0$ for all $m$ then
$$
\liminf_{N\to\infty} \hat{m}(N)\geq m_0.\ \mathbb{P}_{\psi_0}-a.s.
$$
\end{teor}

\noindent{\bf Proof:} From Proposition \ref{douc}-(i) we have:
$$l(\psi_0)-l(\psi)\to H(\psi_0)-H(\psi),$$
$\psi\in\Psi_m$, and $H(\psi)-H(\psi_0)=
\mathbb{E}_{\psi_0}
\left(\log\frac{p_{\psi_0}(y_0|y_{-\infty:-1})}{p_\psi(y_0|y_{-\infty:-1})}\right)
:=D(\psi_0,\psi).
$

Therefore for $m<m_0$:
$$
\inf_{\psi\in\Psi_m}[l(\psi_0)-l(\psi)]\to
l(\psi_0)-l(\hat{\psi})\to
\inf_{\psi\in\Psi_m} D(\psi_0,\psi)>0,
$$
$D(\psi_0,\psi)>0$ since $m_0$ in minimal. We have:
$$\lim_{N\to\infty}l(\hat{\psi}_{m_0})-
l(\hat{\psi}_m)=D(\psi_0,\psi)>0.
$$
By the definition of $C(N,m)$ and by assumption
$\lim_{N}pen(N,m)=0$,
$$\lim_{N\to\infty}C(N,m)-C(N,m_0)=D(\psi_0,\psi)>0,
$$
for any $m<m_0$. On the other hand
$C(N,\hat{m}(N))-C(N,m_0)\leq0$,
by the definition of $\hat{m}(N)$ and  we conclude that
$$\liminf_{N\to\infty} \hat{m}(N)\geq m_0.\ \mathbb{P}_{\psi_0}-a.s
$$.\cqd

In the following we prove that the
estimator PML under-estimate the number state.

Let us define the distribution,
$$
Q_m(y_{1:N}|x_1)=\mathbb{E}_{p(\psi)}(p(y_{1:N}|y_0,x_1,\psi)),
$$
where $p(\psi)$ is a priori distribution on  $\Psi_m$.
In the following we will write the model in its vectorial form,
\begin{equation}
\boldsymbol{y}=\boldsymbol{Z}\theta+\boldsymbol{\varepsilon},
\label{modelovectorial}
\end{equation}
where
$\boldsymbol{\varepsilon}=(\sigma\varepsilon_1,
\ldots,\sigma\varepsilon_N)$, $\boldsymbol{y}=y_{1:N}^t$,
in the case AR-MR $\theta=((\theta_{0,1},\theta_{1,1}),
\ldots,(\theta_{0,m},\theta_{1,m}))^t$,
$$
 \boldsymbol{Z}=\left(\begin{array}{ccc}
  (1,y_0)\car_1(x_1)&
  \cdots& (1,y_0)\car_m(x_1)\\
   \vdots &   & \vdots \\
    (1,
    y_{N-1})
    \car_1(x_N)&\cdots&(1,y_{N-1})
    \car_m(x_N)
  \end{array} \right),
$$
while in the case HMMs $\theta=(\theta_{1},
\ldots,\theta_{m})^t$
$$
 \boldsymbol{Z}=\left(\begin{array}{ccc}
  \car_1(x_1)&
  \cdots& \car_m(x_1)\\
   \vdots &   & \vdots \\
    \car_1(x_N)&\cdots&
    \car_m(x_N)
  \end{array} \right).
$$
Given $x_1,y_0$, the likelihood function
for the model
(\ref{modelovectorial}) is,
    \begin{eqnarray}
    \nonumber
        {p}(\boldsymbol{y}|x_1,y_0,\psi)  &=& \sum_{x_{2:N}\in
        \{1,\ldots,m\}^{N-1}}{p}(\boldsymbol{y},X_{2:N}=x_{2:N}|x_1,\psi)\\
    \label{verovectorial}
        &=& \sum_{x_{2:N}\in
        \{1,\ldots,m\}^{N-1}}{p}(\boldsymbol{y}|X_{2:N}=x_{2:N},x_1,\psi){p}(X_{2:N}=x_{2:N}|x_1,\psi)
    \end{eqnarray}
with,
    \begin{eqnarray*}
        {p}(\boldsymbol{y}|X_{2:N}=x_{2:n},x_1,\psi) &=&
        \mathcal{N}(\boldsymbol{y}-Z\theta|0,\sigma^2I_N)\\
        {p}(X_{2:N}=x_{2:N}|x_1,\psi) &=&
        a_{x_{1}x_{2}}\ldots a_{x_{N-1}x_{N}}.\\
        \end{eqnarray*}
Suppose the following structure of  dependence for
the components $\psi$,
    $$p(\psi)=\left(\prod_{i\in E}p(A_i)\right)
     p(\theta|\sigma^2)p(\sigma^2)$$
and suppose  the following densities that are priors  conjugated
for likelihood function (\ref{verovectorial}):
\begin{enumerate}
    \item{$$\theta\sim\mathcal{ N}(\theta|0,\sigma^2\Sigma)
     =(2\pi\sigma^2)^{-m/2}\det(\Sigma)^{-1/2}\exp\left(-\frac{1}{2\sigma^2}
    \theta^t\Sigma^{-1} \theta\right)$$}
    \item{For $\sigma^2$ is proposed an inverted gamma $\mathcal{IG}$,
    $$\sigma^2\sim
    \mathcal{IG}(v_0/2,u_0/2)=\frac{u_0^{v_0/2}}{{2^{v_0/2}}\Gamma(v_0/2)}\left({\sigma^2}\right)^{-(v_0/2+1)}
  \exp\left(-\frac{u_0}{2\sigma^2}\right),$$
     $\Gamma(u)=\int_0^\infty s^u e^{-s}ds$.}
    \item{$A_i\sim \mathcal{D}(e)$. $\mathcal{D}$
    denotes a Dirichlet density with parameter vector
    $e=(1/2,\ldots,1/2)$,
        $$\mathcal{D}(e)=\frac{\Gamma(m/2)}{
        \Gamma(1/2)^m}\prod_{j=1}^ma_{ij}^{-1/2}.$$}
\end{enumerate}

The following Lemma gives a  bound of
the  likelihood function normalized by  $Q_m$.

\begin{lema}
The prior distribution
$p(\psi)$ satisfies for all $m$ and all $\boldsymbol{y}\in\mathbb{R}^N$
the following  inequalities,
    \begin{eqnarray*}
    \lefteqn{\log\frac{p(\boldsymbol{y}|y_0,x_1,\psi)}{Q_m(\boldsymbol{y}|x_1)}\leq\frac{m(m-1)}{2}\log(N)+c_m(N)+
    \log\Gamma(\frac{u_0}{2})}\\
    &+&\frac{\log\det(\Sigma)}{2}
        +\frac{(1+v_0)}{2}\log(u_0+
    \boldsymbol{y}^tP\boldsymbol{y})
    -\frac{N}{2}-\frac{\log\det(M)}{2}-\log\Gamma\left(\frac{N+v_0}{2}\right)
    \end{eqnarray*}
where $M^{-1}=Z^tZ+\Sigma^{-1}$, $P=I-ZM^tZ^t$ and for $N\geq 4$,
    $$c_m(N)=-m\left(\log\frac{\Gamma(m/2)}{\Gamma(1/2)}-
    \frac{m(m-1)}{4N}+\frac{1}{12N}\right).$$
\end{lema}

Lemma 1 constitutes a  basic step  in the proof of the following
proposition,
\begin{prop}
Let $\hat{m}$ the PML number state. Then for all
$m_0$, all $\psi\in\Psi_{m_0}$ and all $m>m_0$:
$$\mathbb{P}(\hat{m}>m_0)
    \leq\sum_{m=m_0+1}^{m_{\max}}\exp(I'+\Delta pen(m_0,m))
    \int_{\{\boldsymbol{y}\}}(u_0+
    \boldsymbol{y}^tP\boldsymbol{y})^{\frac{N+v_0}{2}}Q_m(\boldsymbol{y}|x_1)
    d\boldsymbol{y}$$
where $\Delta_N pen(m_1,m_2):=pen(N,m_1)-pen(N,m_2)$,
\begin{eqnarray*}
I'&=&\frac{m(m-1)}{2}\log(N)+c_m(N)+
    \log\Gamma(\frac{u_0}{2})\\
    &+&\frac{\log\det(\Sigma)}{2}
    -\frac{N}{2}-\frac{\log\det(M)}{2}-\log\Gamma\left(\frac{N+v_0}{2}\right)
\end{eqnarray*}

\end{prop}
\noindent {\bf Proof:} by Lemma 1,
$$\log\frac{p(\boldsymbol{y}|y_0,x_1,\psi)}{Q_m(\boldsymbol{y}|x_1)}\leq I $$
with
\begin{eqnarray*}
I&=&\frac{m(m-1)}{2}\log(N)+c_m(N)+
    \log\Gamma(\frac{u_0}{2})\\
    &+&\frac{\log\det(\Sigma)}{2}
        +\frac{(1+v_0)}{2}\log(u_0+
    \boldsymbol{y}^tP\boldsymbol{y})
    -\frac{N}{2}-\frac{\log\det(M)}{2}-\log\Gamma\left(\frac{N+v_0}{2}\right),
\end{eqnarray*}
also,
$$\mathbb{P}(\hat{m}(N)>m_0)=\sum_{m=m_0+1}^{m_{\max}}\mathbb{P}(\hat{m}(N)=m),$$
and therefore,
\begin{eqnarray*}
\lefteqn{\mathbb{P}(\hat{m}(N)=m)}\\
&\leq&\mathbb{P}\left(\log p(\boldsymbol{y}|y_0,x_1,\psi_0)
\leq\sup_{\psi\in\mathcal{M}}\log p(\boldsymbol{y}|y_0,x_1,\psi)+\Delta pen(m_0,m)\right)\\
&\leq&\mathbb{P}\left(\log p(\boldsymbol{y}|y_0,x_1,\psi_0)
\leq\log Q_m(\boldsymbol{y}|x_1)+I+\Delta pen(m_0,m)\right)\\
&=&\mathbb{E}(\car_{\{\log p(\boldsymbol{y}|y_0,x_1,\psi_0)
\leq\log Q_m(\boldsymbol{y}|x_1)+I+\Delta pen(m_0,m)\}})\\
&=&\int_{\{\boldsymbol{y}\}}\car_{\{\log p(\boldsymbol{y}|y_0,x_1,\psi_0)
\leq\log Q_m(\boldsymbol{y}|x_1)+I+\Delta pen(m_0,m)\}}(\boldsymbol{y})\exp\log
p(\boldsymbol{y}|y_0,x_1,\psi)d\boldsymbol{y}\\
&\leq&\int_{\{\boldsymbol{y}\}}\exp(\log Q_m(\boldsymbol{y}|x_1)+I+\Delta pen(m_0,m))d\boldsymbol{y}.
\end{eqnarray*}
 get:
    $$\mathbb{P}(\hat{m}>m_0)
    \leq\sum_{m=m_0+1}^{m_{\max}} \exp(I'+\Delta pen(m_0,m))
    \int_{\{\boldsymbol{y}\}}(u_0+
    \boldsymbol{y}^tP\boldsymbol{y})^{\frac{1+v_0}{2}}Q_m(\boldsymbol{y}|x_1)
    d\boldsymbol{y}.$$\cqd \\
As a consequence of this result and the first Borel-Cantelli
Lemma,
the convergence of $\hat{m}$ depends on the study of the series
$\sum_N\sum_{m=m_0+1}^{m_{\max}}\exp(I'(N,m)+\Delta pen(m_0,m))$.

In the following theorem we find under-estimate
 estimator of number state $m_0$.
\begin{teor}
If  $\int_{\{\boldsymbol{y}\}}(u_0+
    \boldsymbol{y}^tP\boldsymbol{y})^{\frac{N+v_0}{2}}Q_m(\boldsymbol{y}|x_1)
    d\boldsymbol{y}<\infty$ and
$
    \lim_{N\to\infty} pen(N,m)-pen(N+1,m)=0,
$
then
$$
\hat{m}(N)\leq m_0\ \ c.s-\mathbb{P}_{\psi_0}.
$$
\end{teor}
\noindent{\bf Proof:} Let us  defined
 $a_N=I'(N,m)+\Delta pen(m_0,m)$. Observe that
the serie
    $$\sum_{m=m_0+1}^M\sum_N\exp a_N<\infty,$$
converges as consequence of the ratio criterio and this
shows that $lim_{N\to\infty}a_{N+1}-a_{N}<0$. In fact,
$$\lim_{N\to\infty}\frac{m(m-1)}{2}\log\left(\frac{N+1}{N}\right)
+c_m(N+1)-c_m(N)=0,$$
$-\frac{1+\log2}{2}<0$,
$-\log\left(\frac{\Gamma((N+1+v_0)/2)}{\Gamma((N+v_0)/2)}\right)<0$

\begin{eqnarray*}
\lefteqn{\lim_{N\to\infty}\Delta_{N+1} pen(m_0,m)+\Delta_N pen(m_0,m)}\\
&=&\lim_{N\to\infty} pen(N+1,m_0)-pen(N+1,m)-pen(N,m_0)+pen(N,m)\\
&\leq&\lim_{N\to\infty} pen(N,m)-pen(N+1,m)=0.
\end{eqnarray*}
Then we have
\begin{eqnarray*}
\lefteqn{\lim_{N\to\infty}a_{N+1}-a_{N}}\\
&=&\lim_{N\to\infty}\frac{m(m-1)}{2}\log\left(\frac{N+1}{N}\right)
+c_m(N+1)-c_m(N)-\frac{1+\log2}{2}\\
&-&
\log\left(\frac{\Gamma((N+1+v_0)/2)}{\Gamma((N+v_0)/2)}\right)
+\Delta_{N+1} pen(m_0,m)-\Delta_N pen(m_0,m)<0
\end{eqnarray*}

Thus $\sum_N\mathbb{P}_{\psi_0}(\hat{m}(N)>m_0)<\infty$ and from
the Borel-Cantelli
lemma we conclude  that $\mathbb{P}_{\psi_0}(\hat{m}(N)>m_0\ i.o)=0$.
This is  equivalent to say that $\hat{m}(N)\leq m_0$ c.s-$\mathbb{P}_{\psi_0}$.
\cqd\\

One of the most common choices is $pen(N,m)=\frac{\log(N)}{2}dim(\Psi_m)$
(Bayesian information criteria, BIC). It is natural to use
$dim(\Psi_m)=m(m-1)+m dim(\Theta)+1$.\\

\appendix

\section{Proof of Lemma 1}
The proof of this Lemma is obtained by showing the existence of  constants
$C_1,C_2$ such that:
\begin{align}
p(\boldsymbol{y}|x_1,\psi,)&\leq C_1Q_m(\boldsymbol{y}|x_{1:N})\\
p(x_{2:N}|x_1,\psi)&\leq C_2Q_m(x_{2:N}|x_1).
\end{align}
This would implies that,
\begin{eqnarray*}
{p}(\boldsymbol{y}|y_0,x_1,\psi)
&=& \sum_{x\in
    \{1,\ldots,m\}^N}{p}(\boldsymbol{y}|x_{1:N},\psi){p}(x_{2:N}|x_1,\psi)\\
&\leq& C_1C_2 \sum_{x\in
    \{1,\ldots,m\}^N}Q_m(\boldsymbol{y}|x_{2:N})Q_m(x_{2:N}|x_1)\\
&=&C_1C_2Q_m(\boldsymbol{y}|x_1).
\end{eqnarray*}
and hence ${p}(\boldsymbol{y}|y_0,x_1,\psi)\leq C_1C_2Q_m(\boldsymbol{y}|x_1)$.

We proceed with the evaluation of
$Q_m(x_{2:N}|x_1)$ following the proof given in the appendix of \cite{Liu}.
Let
$$Q_m(x_{2:N}|x_1)=\prod_{i=1}^m\left[\frac{\Gamma(m/2)}{\Gamma(N_i+1/2)}
\left(\prod_{i=1}^m\frac{\Gamma(N_{ij}+1/2)}{\Gamma(1/2)}\right)\right]$$
and
   \begin{equation}
     \frac{p(x_{2:N}|x_1,\psi)}{Q_m(x_{2:N}|x_1)}
     \leq\frac{\prod_{i=1}^m\prod_{j=1}^m{(\frac{N_{ij}}{N_i})}^{N_{ij}}}{\prod_{i=1}^m
     \left[\frac{\Gamma(m/2)}{\Gamma(N_i+1/2)}
     \left(\prod_{i=1}^m\frac{\Gamma(N_{ij}+1/2)}{\Gamma(1/2)}\right)\right]}.
     \label{cotax1}
    \end{equation}
We have that
and the right side the equation (\ref{cotax1}) does not exceed,
    $$\left[\frac{\Gamma(N+m/2)\Gamma(1/2)}{\Gamma(m/2)\Gamma(N+1/2)}
    \right]^m.$$
In Gassiat and Boucheron \cite{gassiat3},
is noted that,
    $$m\log\left[\frac{\Gamma(N+m/2)\Gamma(1/2)}{\Gamma(m/2)\Gamma(N+1/2)}
    \right]\leq\frac{m(m-1)}{2}\log N+c_m(N),$$
for $N\geq4$, $c_m(N)$ is choosed as:
    $$-m\left(\log\frac{\Gamma(m/2)}{\Gamma(1/2)}-
    \frac{m(m-1)}{4N}+\frac{1}{12N}\right).$$
Then:
    \begin{equation}
    \frac{p(\boldsymbol{x}|x_1,\psi)}{Q_m(\boldsymbol{x}|x_1)}
        \leq N^{m(m-1)/2}\exp{c_m(N)}.
    \label{cotax2}
    \end{equation}
To evaluate $Q(\boldsymbol{y}|x_{1:N})$ let us  develop the expression,
\begin{eqnarray*}
p(\boldsymbol{y}|y_0,x_{1:N},\theta,\sigma^2)p(\theta|\sigma^2)p(\sigma^2)
&=&\mathcal{N}(\boldsymbol{y}-Z\theta|0,\sigma^2I_N)\mathcal{ N}(\theta|0,\sigma^2\Sigma)
\mathcal{IG}(\sigma^2|v_0/2,u_0/2)\\
&=&(2\pi\sigma^2)^{-N/2}\exp\left(-\frac{1}{2\sigma^2}
    (\boldsymbol{y}-Z\theta)^t(\boldsymbol{y}-Z\theta)\right)\\
&&(2\pi\sigma^2)^{-m/2}\det(\Sigma)^{-1/2}\exp\left(-\frac{1}{2\sigma^2}
    \theta^t\Sigma^{-1}\theta\right)\\
& &\frac{u_0^{v_0/2}}{{2^{v_0/2}}\Gamma(v_0/2)}\left({\sigma^2}\right)^{-(v_0/2+1)}
  \exp\left(-\frac{u_0}{2\sigma^2}\right)
\end{eqnarray*}
The above-mentioned is equivalent to
$$
    \frac{u_0^{v_0/2}(2\pi\sigma^2)^{-N/2}(2\pi\sigma^2)^{-m/2}}
    {{2^{v_0/2}}\Gamma(v_0/2)}\exp\left(-\frac{(\theta-\boldsymbol{m})^tM^{-1}
    (\theta-\boldsymbol{m})}{2\sigma^2}
    \right)
    \left({\sigma^2}\right)^{-(v_0/2+1)}\exp\left(-\frac{(u_0+\boldsymbol{y}^tP\boldsymbol{y})}
    {2\sigma^2}\right)
 $$
with $M^{-1}=Z^tZ+\Sigma^{-1}$, $\boldsymbol{m}=MZ^t\boldsymbol{y}$ and $P=I-ZM^tZ^t$.
Integrating the last expression respect to $\theta$ and then
to $\sigma^2$ we obtain

    \begin{equation}
    Q(\boldsymbol{y}|x_{1:N})=\frac{u_0^{v_0/2}\det(M)^{1/2}\Gamma((N+v_0)/2)}
    {(\pi\sigma^2)^{N/2}\Gamma(v_0/2)\det(\Sigma)^{1/2}(u_0+
    \boldsymbol{y}^tP\boldsymbol{y})^{(N+v_0)/2}},
    \end{equation}
this given,
    \begin{eqnarray*}
    \frac{p(\boldsymbol{y}|y_0,x_1,\psi)}{Q_m(\boldsymbol{y}|x_{1:N})}
    &\leq&\frac{p(\boldsymbol{y}|y_0,x_1,\hat{\psi})(\pi\sigma^2)^{N/2}\Gamma(v_0/2)\det(\Sigma)^{1/2}(u_0+
    \boldsymbol{y}^tP\boldsymbol{y})^{(N+v_0)/2}}{u_0^{v_0/2}\det(M)^{1/2}\Gamma((N+v_0)/2)}\\
    &=&\frac{\exp\left(-\frac{1}{2\hat{\sigma^2}}
    (\boldsymbol{y}-Z\hat{\theta})^t(\boldsymbol{y}-Z\hat{\theta})\right)
    (\pi\hat{\sigma^2})^{N/2}\Gamma(v_0/2)\det(\Sigma)^{1/2}(u_0+
    \boldsymbol{y}^tP\boldsymbol{y})^{(N+v_0)/2}}{(2\pi\hat{\sigma^2})^{N/2}
    u_0^{v_0/2}\det(M)^{1/2}\Gamma((N+v_0)/2)}\\
    &=&\frac{\exp\left(-\frac{N-1}{2}
    \right)\Gamma(v_0/2)\det(\Sigma)^{1/2}(u_0+
    \boldsymbol{y}^tP\boldsymbol{y})^{(N+v_0)/2}}{2^{N/2}u_0^{v_0/2}\det(M)^{1/2}\Gamma((N+v_0)/2)}.
    \end{eqnarray*}
with this expression and the equation  (\ref{cotax2}) we obtain lemma 1.\cqd

\section{Simulations}
In this section we apply our results to some simulated data. We work with an HMMs and
two AR-RM. We use  $pen=\frac{\log(N)}{2}dim(\Psi_m)$ (BIC).
We value the likelihood function for any set of parameters $\psi$ by computing
$$
  p(y_{1:N}|y_0\psi)=\sum_{i=1}^m\alpha_N(i),
$$
where $\alpha_n(i)=p(y_{1:n},X_n=i)$ can be evaluated
recursively with the following formulae forward of Baum,
$$
\alpha_n(j)=\sum_{i=1}^m\alpha_{n-1}(i)a_{ij}p(y_n|y_{n-1},X_n=i)
$$
see D. Le Nhu {\it et alii} \cite{Nlputerman}.\\

\subsection{HMMs}
In the simulation of the HMMs we set the following parameters:
$dim(\Psi_m)=m^2+1$ $N=500$, $m=3$,
$\sigma^2=1.5$,  $\theta=(-2,\ 1,\ 4)$,
$$
A=\left(\begin{array}{ccc}
   0\mbox{.}9 & 0\mbox{.}05 & 0\mbox{.}05\\
0\mbox{.}05 &  0\mbox{.}9 &0\mbox{.}05 \\
 0\mbox{.}05 & 0\mbox{.}05 &0\mbox{.}9 \\
 \end{array}\right),
$$
the observed serie is plotted in figure 1.

The table 1 contains the values for the penalized maximum
likelihood for $m=2,\ldots,7$, we observe that $\hat{m}=3$. In
this case $\hat{\psi}$ is estimated by using the SAEM, the values
are, $\hat{\sigma^2}=1.49$, $\hat{\theta}=\left(-1.98,\ 4.09,\
0.91\right)$,
$$
\hat{A}=\left(\begin{array}{ccc}
   0.8650 & 0.0274& 0.1076\\
   0.0404 &0.8943 & 0.0653\\
   0.0658&0.0648& 0.8694 \
 \end{array}\right),
$$
in the figure 2 displayed the
sequence $\{\psi^{(t)}\}$, $t=1,\ldots,4000$ and we observe the convergence
of the estimate.

\begin{table}
\begin{center}
\begin{tabular}{|c|c|c|c|}
\hline
$m$ & $-l(\psi)$ & $pen$ & $-l(\psi)+pen$\\ \hline
2 & 802.32 & 15.53 &817.85 \\ \hline
3 & 419.09 &31.07 &450.16 \\ \hline
4 & 417.70 &52.82 &470.52\\ \hline
5 & 464.70 & 80.78 & 545.48\\ \hline
6 & 445.89 &114.97& 560.86\\ \hline
7 & 436.26 & 155.36& 591.62 \\ \hline
\end{tabular}
\end{center}
\caption{The values for the PML}
\end{table}

\begin{figure}
\begin{center}
\includegraphics[angle=0,scale=0.60]{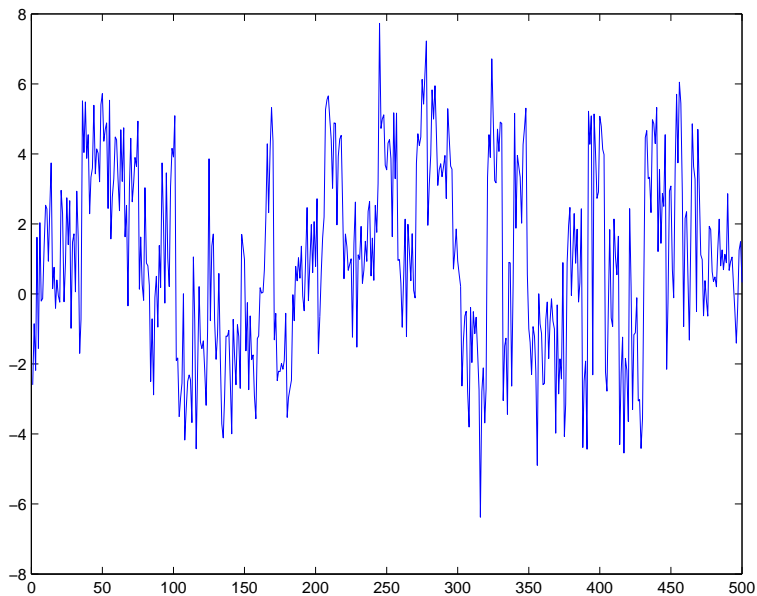}
\end{center}
\caption{The observed serie $y_1,\ldots,y_{500}$ for the HMMs}.
\end{figure}

\begin{figure}
\begin{center}
\begin{tabular}{ccc}
\includegraphics[angle=0,scale=0.32]{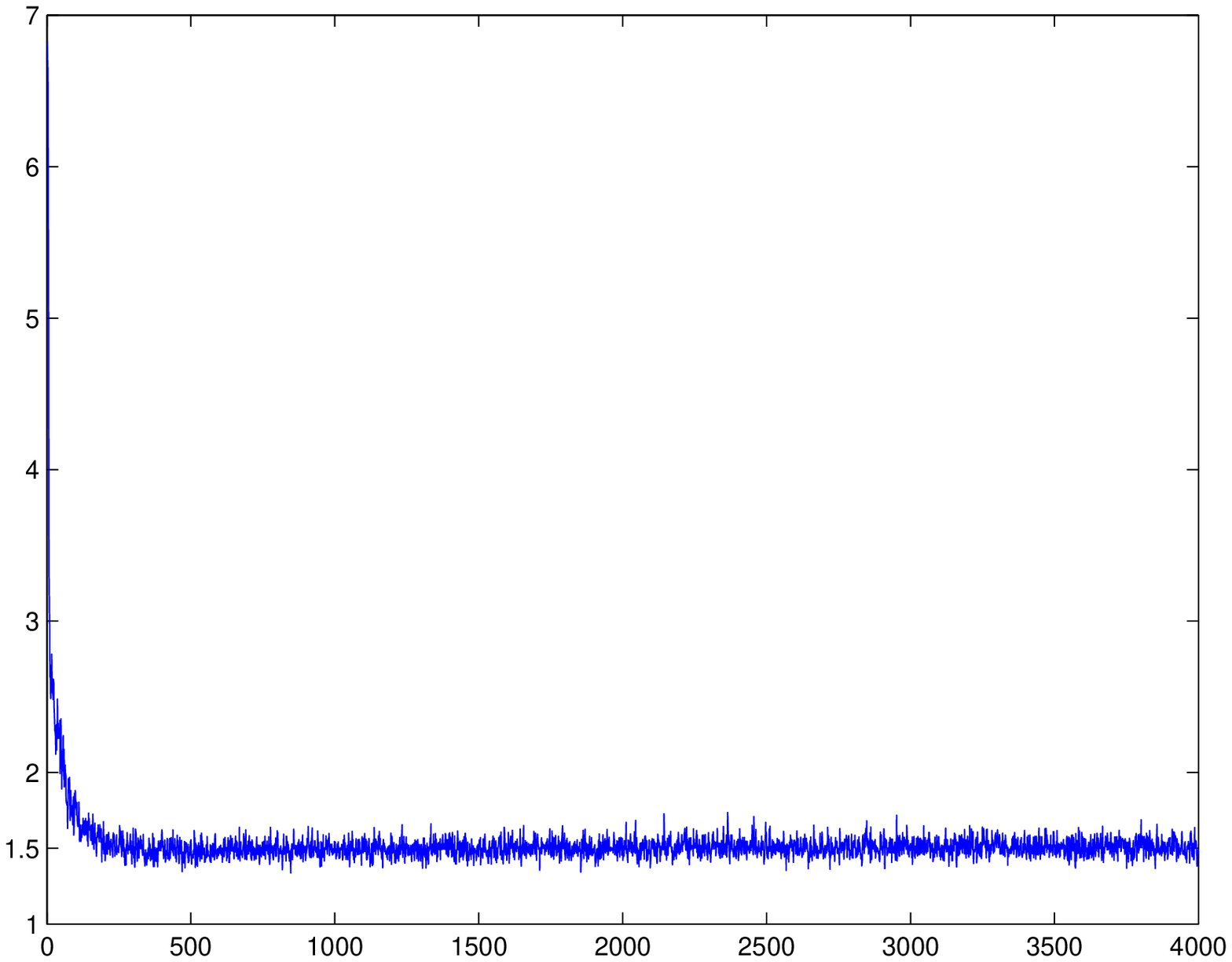}
&
\includegraphics[angle=0,scale=0.32]{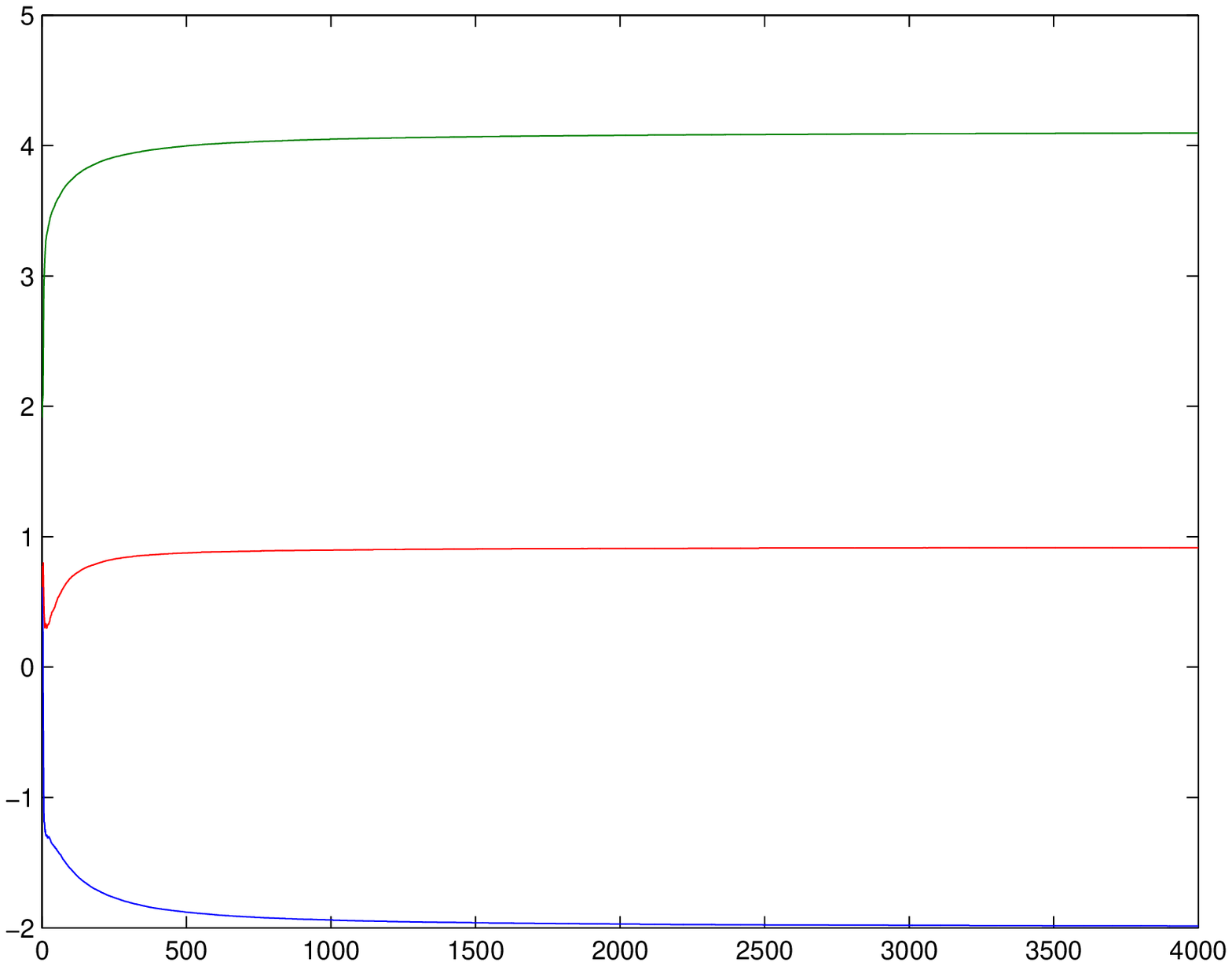}
&
\includegraphics[angle=0,scale=0.32]{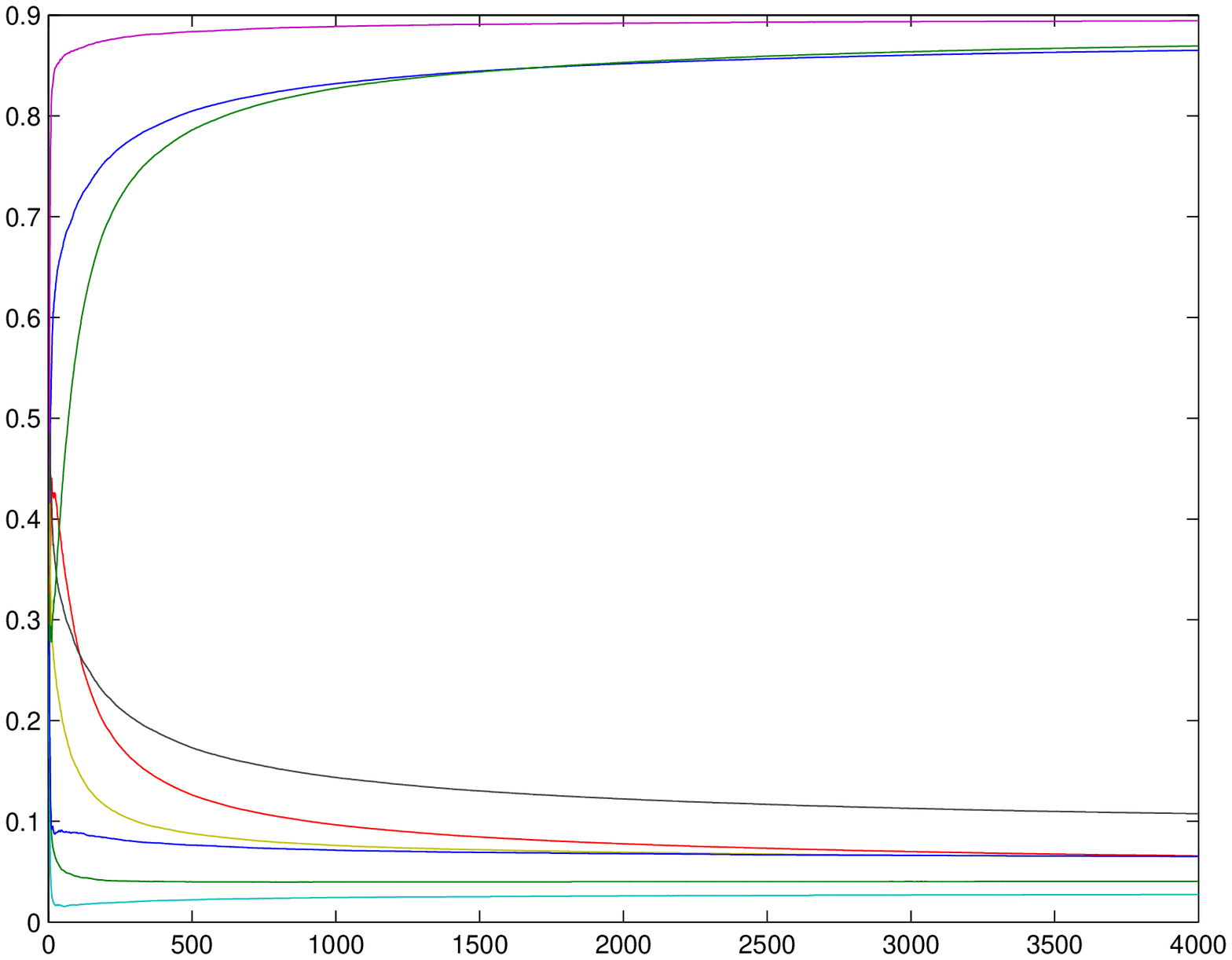}\
\end{tabular}
\end{center}
\caption{Convergence of the estimate of, $\sigma^2$, $\theta$ and $A$.}
\end{figure}

\subsection{AR-RM}
In the first simulation of the AR-RM we set the following parameters: $dim(\Psi_m)=m(m+1)+1$,
$N=500$, $m=2$, $\sigma^2=1\mbox{.}5$,
\begin{center}
\begin{tabular}{cc}
 $\theta=\left(\begin{array}{cc}
   1 & -1 \\
   -0\mbox{.}5 &   0\mbox{.}5\
 \end{array}\right),$ &
 $A=\left(\begin{array}{cc}
   0\mbox{.}9 & 0\mbox{.}1 \\
   0\mbox{.}1 & 0\mbox{.}9 \
 \end{array}\right), $
  \\
\end{tabular}
\end{center}
the observed serie is plotted in figure 3.

The table 2 contains the values for the penalized maximum
likelihood for $m=2,\ldots,6$, we observe that $\hat{m}=2$. In
this case $\hat{\psi}$ is estimated by using the SAEM, the values
are, $\hat{\sigma^2}=1.42$,
\begin{center}
\begin{tabular}{cc}
  $\hat{\theta}=\left(\begin{array}{cc}
   1\mbox{.}07 & -0\mbox{.}96 \\
   -0\mbox{.}5 & 0\mbox{.}5 \
 \end{array}\right) $ &
$\hat{A}=\left(\begin{array}{cc}
   0.8650 & 0.1350 \\
   0.1130 &   0.8870\
 \end{array}\right).$  \\
\end{tabular}
\end{center}
in the figure 4 displayed the
sequence $\{\psi^{(t)}\}$, $t=1,\ldots,1000$ and we observe the convergence
of the estimate.

\begin{figure}
\begin{center}
\includegraphics[angle=0,scale=0.60]{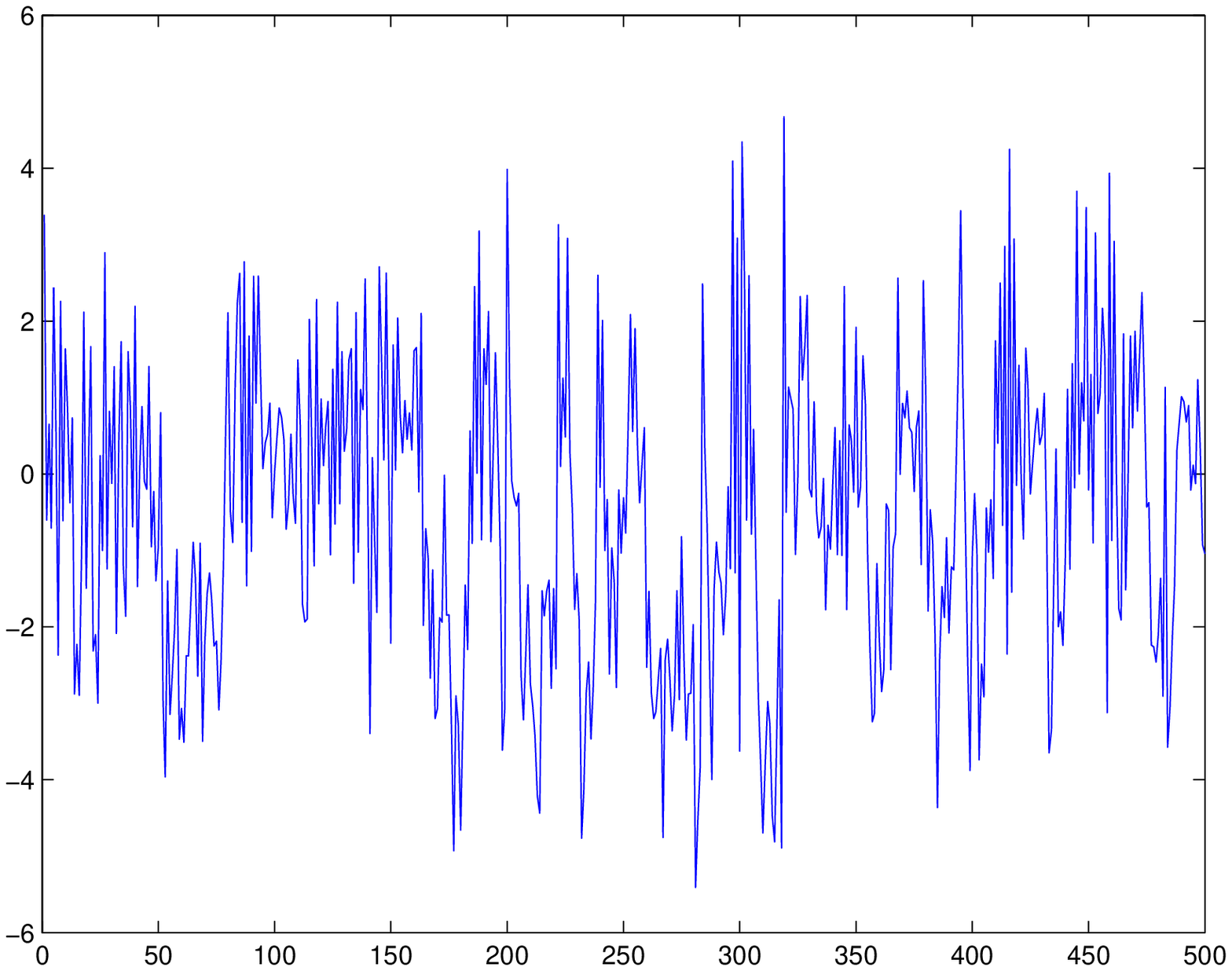}
\end{center}
\caption{The observed serie $y_1,\ldots,y_{500}$ for the AR-MR}
\end{figure}

\begin{table}
\begin{center}
\begin{tabular}{|c|c|c|c|}
\hline
$m$ & $-l(\psi)$ & $pen$ & $-l(\psi)+pen$\\ \hline
2 & 351.14 & 18.64 &369.78 \\ \hline
3 & 346.64 & 37.28 &383.92 \\ \hline
4 & 355.10 &64.14 &417.24\\ \hline
5 & 354.52 & 93.21 & 447.73\\ \hline
6 & 384.50 &130.50& 515.00\\ \hline
\end{tabular}
\end{center}
\caption{The values for the PML}
\end{table}

\begin{figure}
\begin{center}
\begin{tabular}{cc}
\includegraphics[angle=0,scale=0.35]{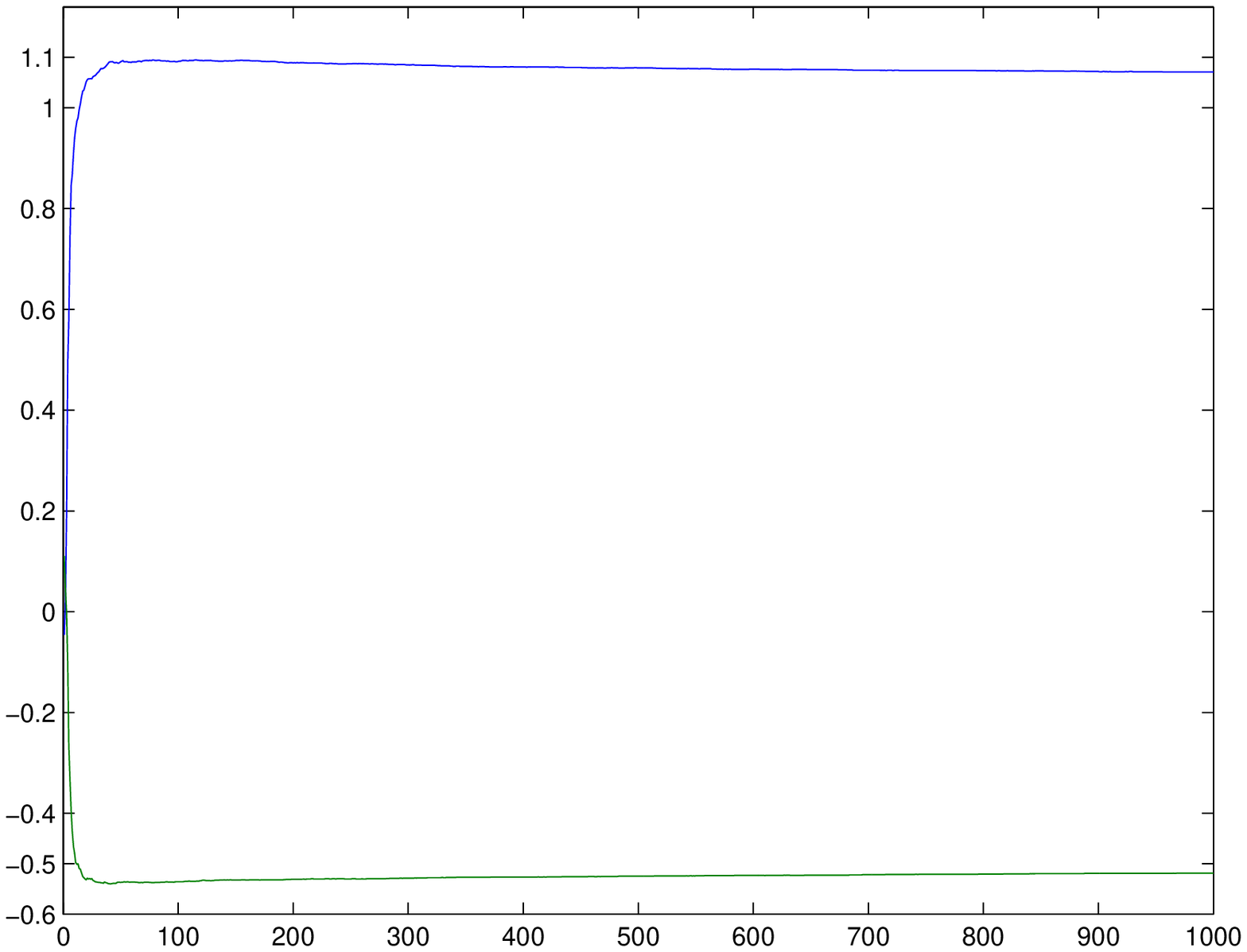}
&
\includegraphics[angle=0,scale=0.35]{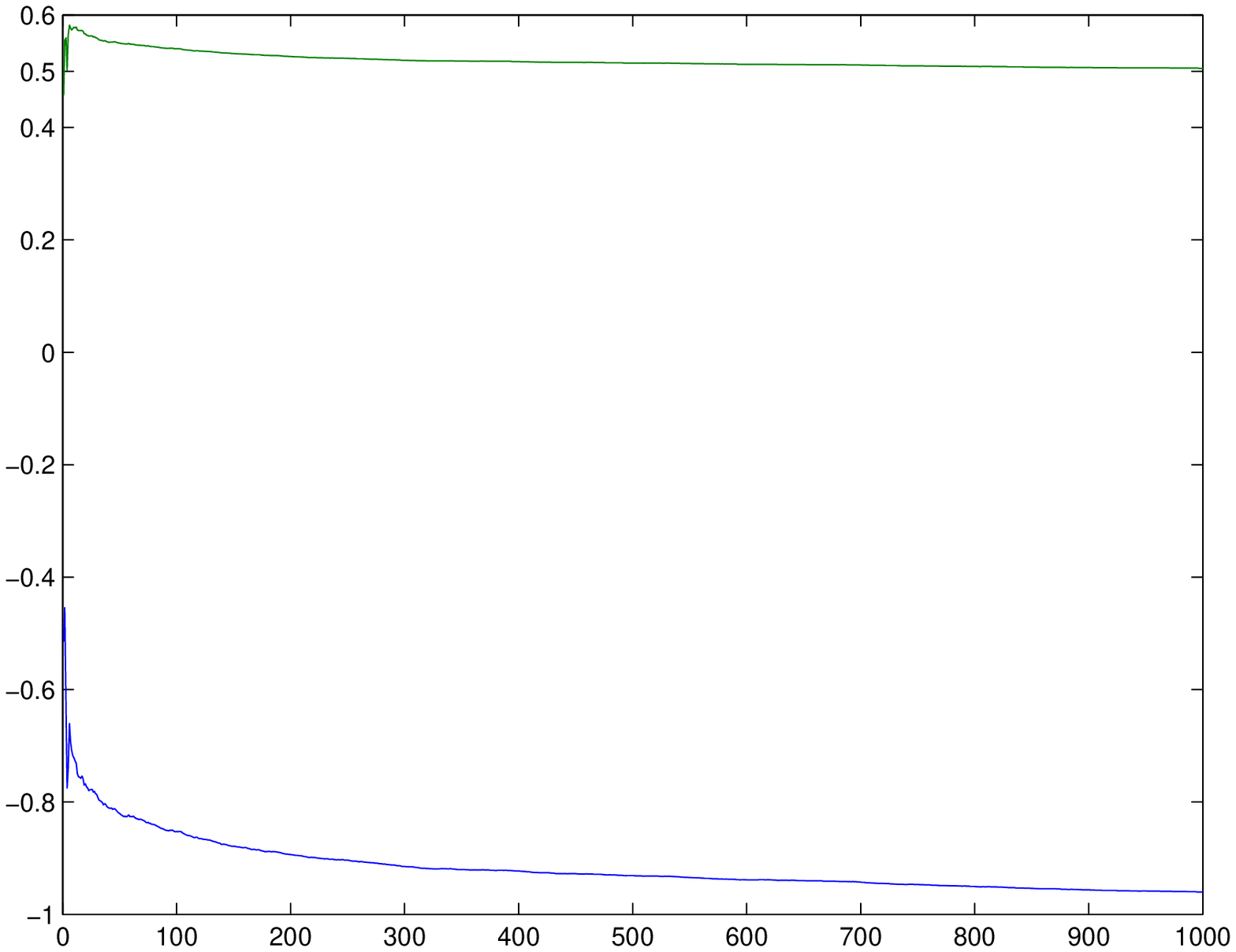} \\
\includegraphics[angle=0,scale=0.35]{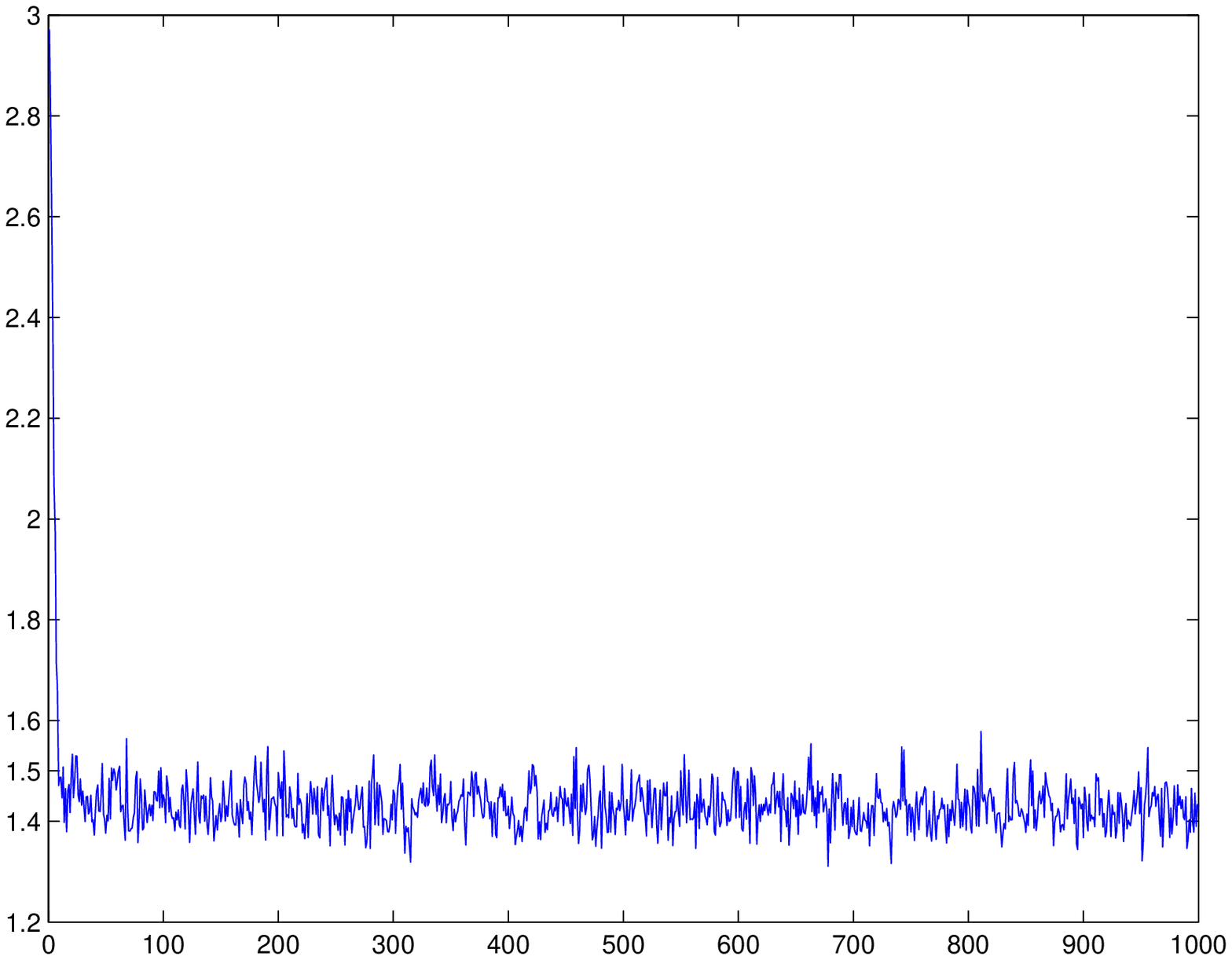}
&\includegraphics[angle=0,scale=0.35]{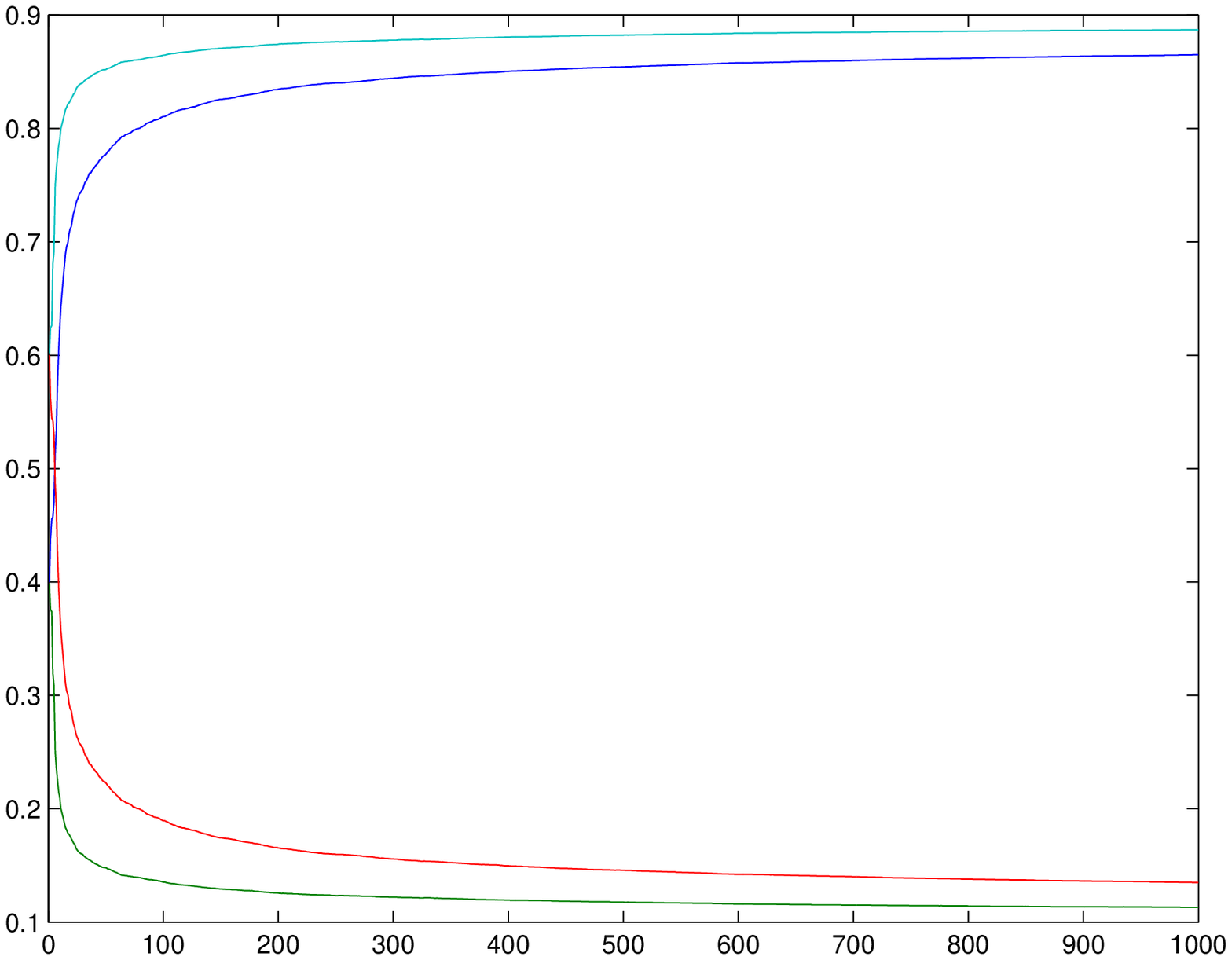}
\end{tabular}
\end{center}
\caption{Convergence of the estimate of,  $\theta_1$, $\theta_2$, $\sigma^2$, and $A$.}
\end{figure}

In the second simulation of the AR-RM we set the following parameters:
$N=500$, $m=2$, $\sigma^2=1\mbox{.}5$,
\begin{center}
\begin{tabular}{cc}
 $\theta=\left(\begin{array}{cc}
   1 & -2 \\
   -0\mbox{.}7 &   1\mbox{.}08\
 \end{array}\right)$ &
 $A=\left(\begin{array}{cc}
   0\mbox{.}9 & 0\mbox{.}1 \\
   0\mbox{.}1 & 0\mbox{.}9 \
 \end{array}\right),$ \\
\end{tabular}
\end{center}
the observed serie is plotted in figure 3.

In this case $m=2$ is fixed and  $\hat{\psi}$ is estimated by using the SAEM, the
values are, $\hat{\sigma^2}=1.42$,
\begin{center}
\begin{tabular}{cc}
  $\hat{\theta}=\left(\begin{array}{cc}
   0\mbox{.}85 & -2\mbox{.}01 \\
   -0.69 & 1.08\
 \end{array}\right) $ &
$\hat{A}=\left(\begin{array}{cc}
   0.9093 & 0.0907 \\
   0.019 &   0.9181\
 \end{array}\right),$  \\
\end{tabular}
\end{center}
in the figure 6 displayed the
sequence $\{\psi^{(t)}\}$, $t=1,\ldots,1000$ and we observe the convergence
of the estimate.\\

\begin{figure}
\begin{center}
\includegraphics[angle=0,scale=0.60]{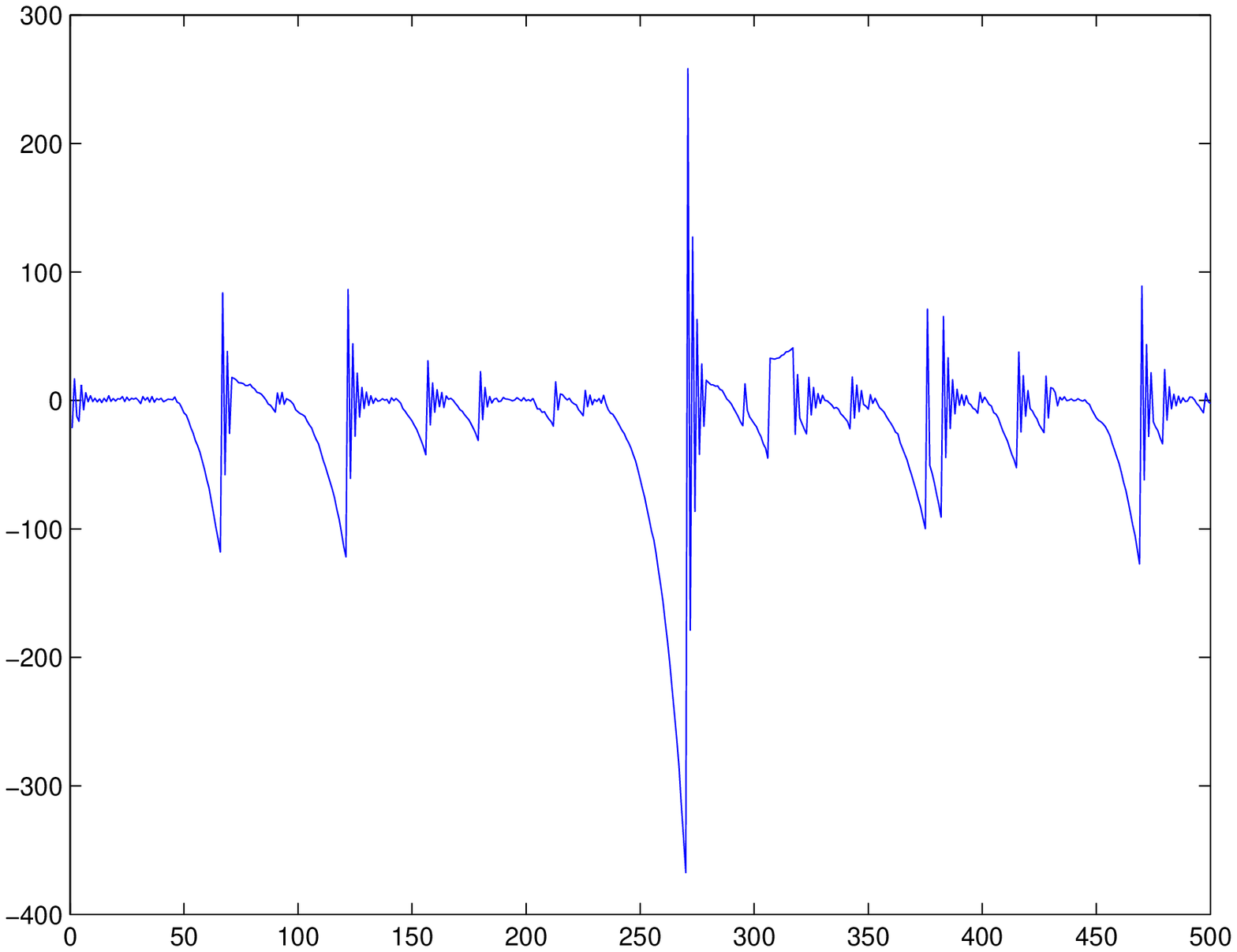}
\end{center}
\caption{The observed serie $y_1,\ldots,y_{500}$ for the AR-MR}.
\end{figure}

\begin{figure}
\begin{center}
\begin{tabular}{cc}
\includegraphics[angle=0,scale=0.35]{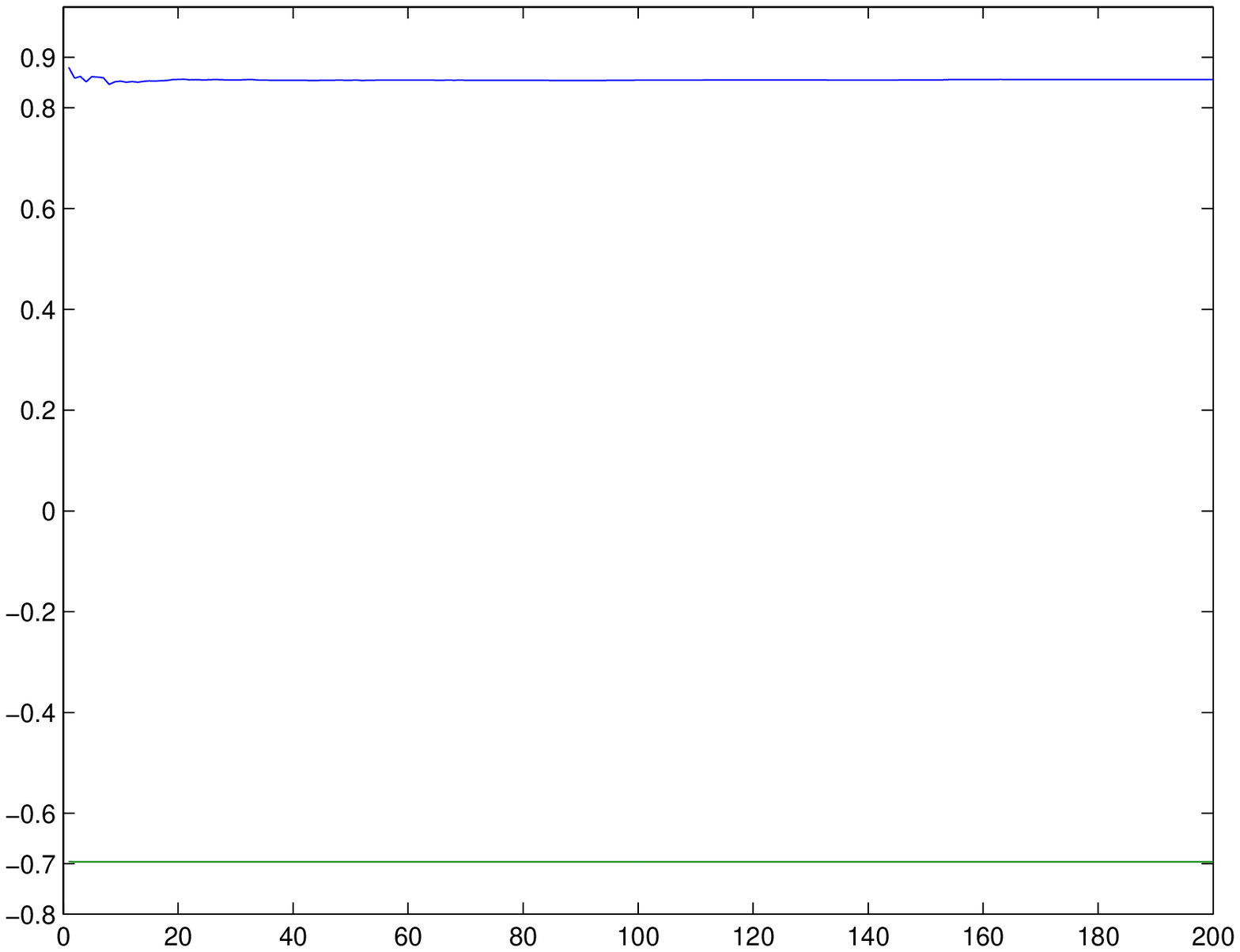}
&
\includegraphics[angle=0,scale=0.35]{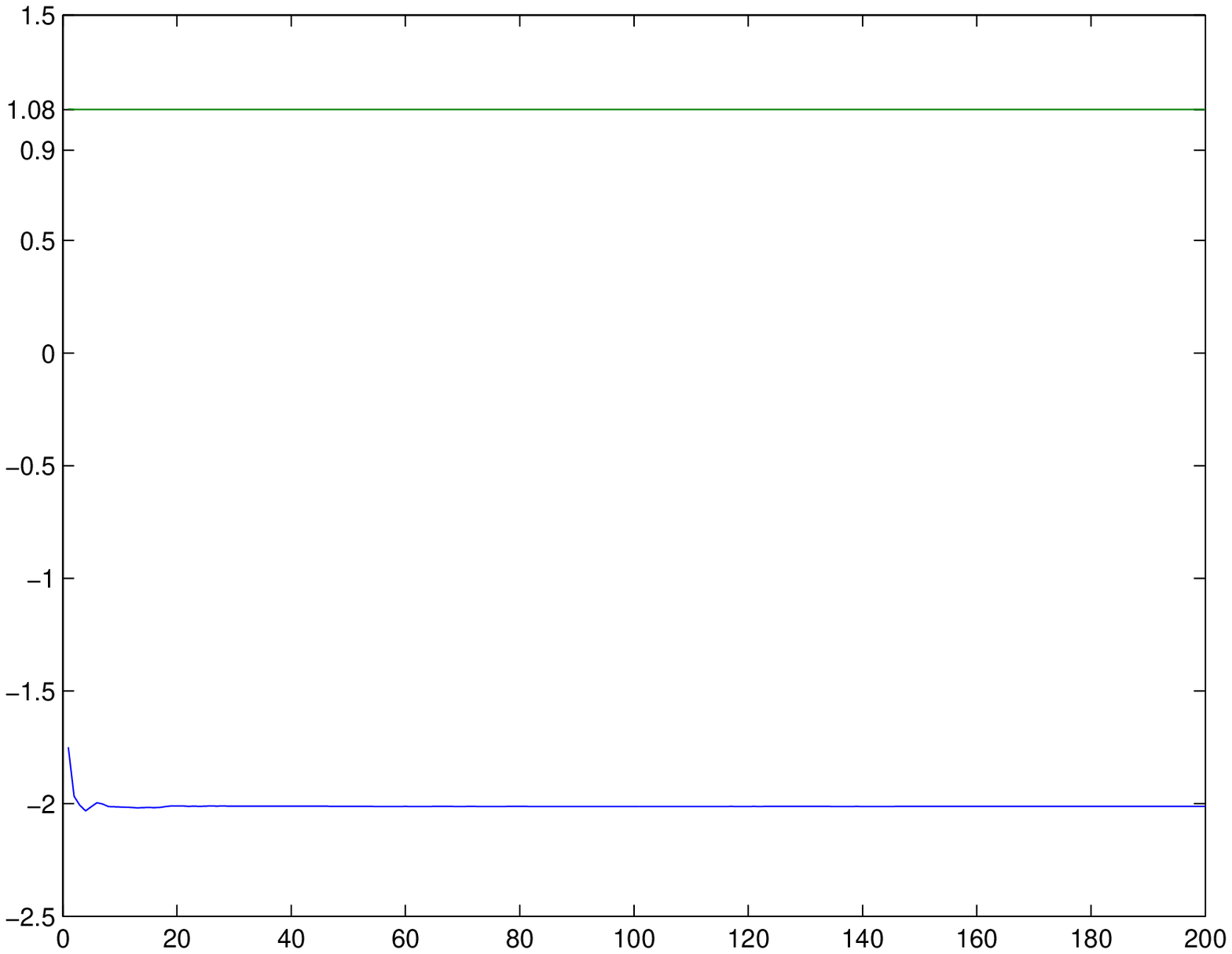} \\
\includegraphics[angle=0,scale=0.35]{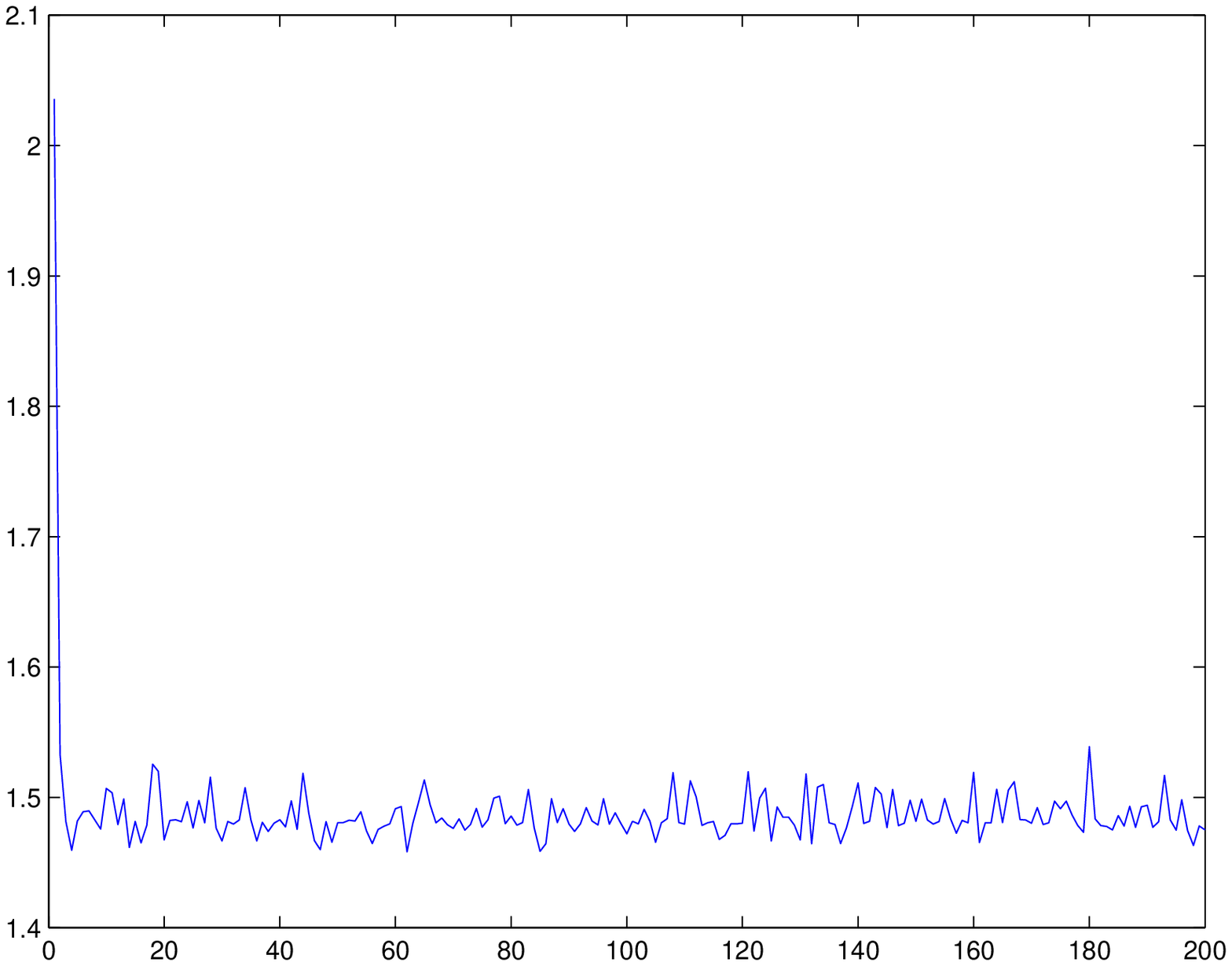}
&\includegraphics[angle=0,scale=0.35]{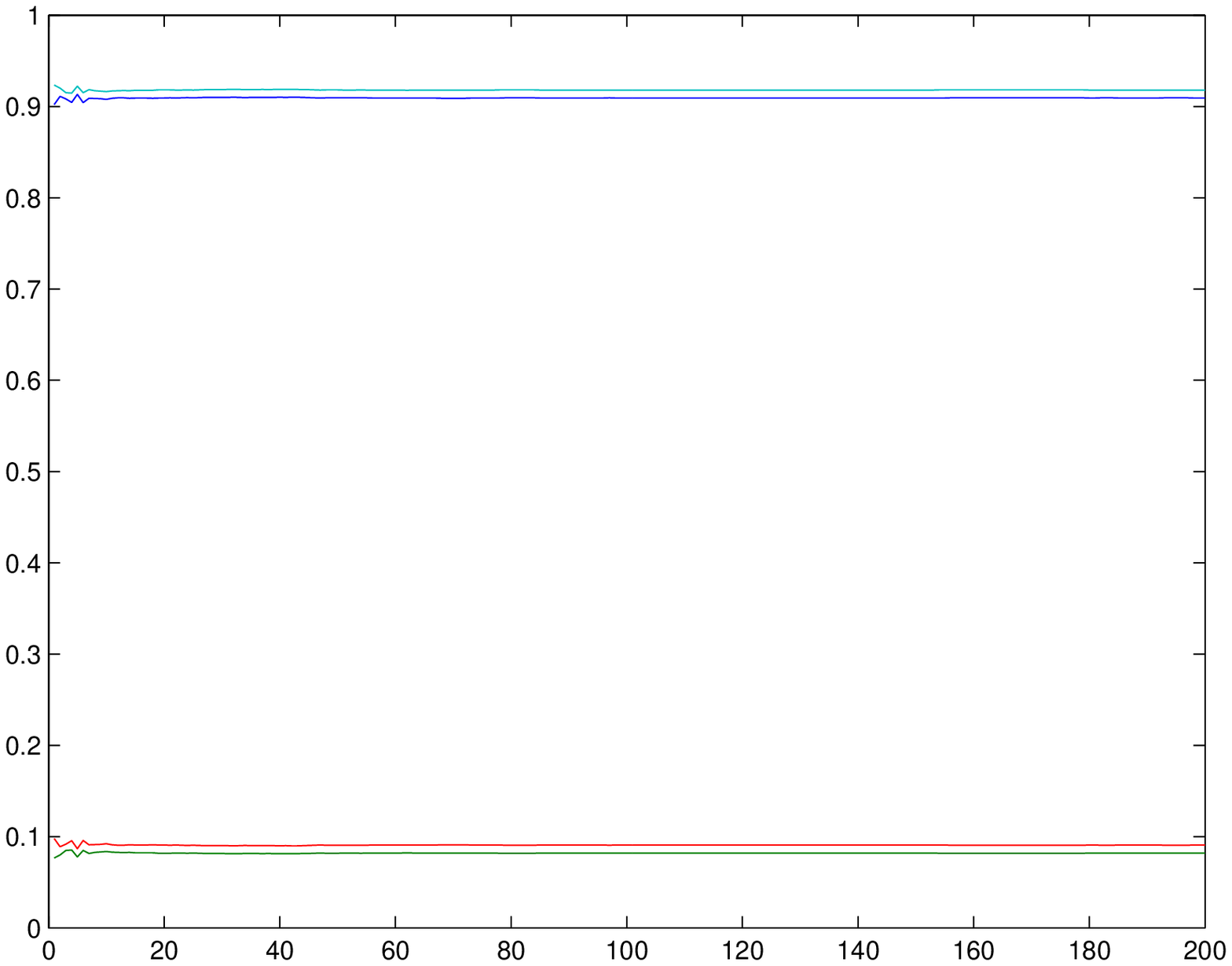}
\end{tabular}
\end{center}
\caption{Convergence of the estimate of,  $\theta_1$, $\theta_2$, $\sigma^2$, and $A$.}
\end{figure}

{\noindent\bf Acknowledgements.} The authors thanks to Marc
Lavielle and Jean-Michel Loubes for help concerning this work. The
authors also thanks Rafael Rosales for carefully reading a
preliminary version.

\bibliographystyle{plain}
\bibliography{arhmm}

\end{document}